\documentclass[12pt]{article}
\usepackage{color}
\usepackage{graphicx} 
\usepackage{amsmath}
\usepackage{amssymb}
\begin{document}
\newcommand{\vare}{\varepsilon} 
\newtheorem{tth}{Theorem}[section]
\newtheorem{dfn}[tth]{Definition}
\newtheorem{lem}[tth]{Lemma}
\newtheorem{prop}[tth]{Proposition}
\newtheorem{coro}[tth]{Corollary}
\renewcommand{\quad}{\hspace*{1em}}
\begin{center}
{\Large {\bf Genericity of Caustics and Wavefronts on an $r$-corner} }
\vspace*{0.4cm}\\
{\large Takaharu Tsukada}\footnote{Higashijujo 3-1-16 
Kita-ku, Tokyo 114-0001
JAPAN. e-mail : tsukada@math.chs.nihon-u.ac.jp}
\vspace*{0.2cm}\\
{\large  College of Humanities \& Sciences, Department of Mathematics,\\
 Nihon University}
\end{center}
\begin{abstract}
We investigate genericities of reticular Lagrangian maps and 
reticular Legendrian maps in order to give generic classifications
of caustics and wavefronts generated by a hypersurface germ
without or with a boundary in a smooth manifold.
We also give simpler proofs of main results in \cite{retLag},\cite{retLeg}.
\end{abstract}
\section{Introduction}
\quad 
Lagrangian and Legendrian  singularities can be found in many problems of differential geometry, 
calculus of variations and mathematical physics. 
One of the most successful their applications are 
the study of singularities of caustics and wavefronts. 
For example, the particles incident along geodesics  
from a smooth hypersurface in a 
Riemannian manifold to conormal directions define 
a Lagrangian submanifold at a point in the cotangent bundle
and define a Legendrian submanifolds at a point in the $1$-jet bundle. 
The caustic generated by the hypersurface is regarded as 
the caustic of the Lagrangian map 
defined by the restriction of the cotangent bundle projection to the 
Lagrangian submanifold 
and the wavefront generated by the hypersurface is regarded as 
the wavefront of the Legendrian map defined by 
the restriction of the $1$-jet bundle projection to the Legendrian submanifold.
Therefore the studies of the caustics and wavefronts 
generated by smooth hypersurfaces are 
reduced to the studies of Lagrangian and Legendrian  singularities. 


In \cite{retLag} and \cite{retLeg} we investigated the more general cases 
when the hypersurface 
has a boundary, a corner, or {\it an $r$-corner}. 
In these cases particles incident from each edge of the hypersurface 
gives a {\it symplectic regular $r$-cubic configuration} 
at a point of the cotangent bundle 
which is a generalisation of the notion of Lagrangian submanifolds and 
particles incident from each edge of the hypersurface 
gives a {\it contact regular $r$-cubic configuration} 
at a point of the $1$-jet bundle
which is a generalisation of the notion of Legendrian submanifolds.
The caustic generated by the hypersurface germ with an $r$-corner is given by 
the caustic of the symplectic regular $r$-cubic configuration 
which is a generalisation of 
the notion of {\it quasicaustics} given by S.Janeszko (cf., \cite{janeszko3}). 
In these papers 
we investigated the stabilities of caustics and wavefronts generated by the hypersurface
with an $r$-corner by studying the stabilities of 
{\it reticular Lagrangian, Legendrian maps} 
which are generalisations of the notions of Lagrangian, Legendrian maps
for our situations.

\vspace{3mm}
In this paper, we investigate the genericities of caustics and wavefronts generated by a hypersurface
with an $r$-corner.
In order to realize this purpose, 
we shall investigate the  genericities of reticular Lagrangian, Legendrian maps.
In these processes, we shall need to prove that the {\it stabilities} and the 
{\it transversally stabilities}
 of reticular Lagrangian, Legendrian maps are equivalent respectively.
These proofs are simpler than that of the assertions (1)$\Leftrightarrow$(5) of 
Theorem 5.5 in \cite[p.587]{retLag} and 
Theorem 7.4 in \cite[p.123]{retLeg} respectively.

The main results in this paper are generic classification of 
caustics generated by the hypersurface with an $r$-corner on an $n$ dimensional 
manifold  in the cases $n\leq 5,r=0$ and $n\leq 3,r=1$
and a generic classification of wavefronts generated by the hypersurface with
an $r$-corner on an $n$ dimensional manifold 
in the cases $n\leq 6,r=0$ and $n\leq 4,r=1$.
In order to realize this, we shall 
classify generic reticular Lagrangian, Legendrian maps 
for the above cases.

By our theory, we have that: 
A generic caustics is one of the types
$A_2,A^\pm_3,A_4,A^\pm_5,A_6,$ $
D^\pm_4,D^\pm_4,D^\pm_6,E^\pm_6$
in the case $n\leq 5,r=0$  and 
$B^\pm_2,B^\pm_3,B^\pm_4,C^\pm_3,C^\pm_4,F^\pm_4$ in 
the case $n\leq 3,r=1$.
A generic wavefront is one of the types
$A_1,A_2,A_3,A_4,A_4,A_6,D^\pm_4,D_5,D^\pm_6,E_6$
in  the case $n\leq 6,r=0$ and 
$B_2,B_3,B_4,C^\pm_3,C_4,F_4$ in the case $n\leq 4,r=1$.

\begin{figure}[htbp]
\[
\begin{array}{cc}
\includegraphics*[width=7cm,height=7cm]{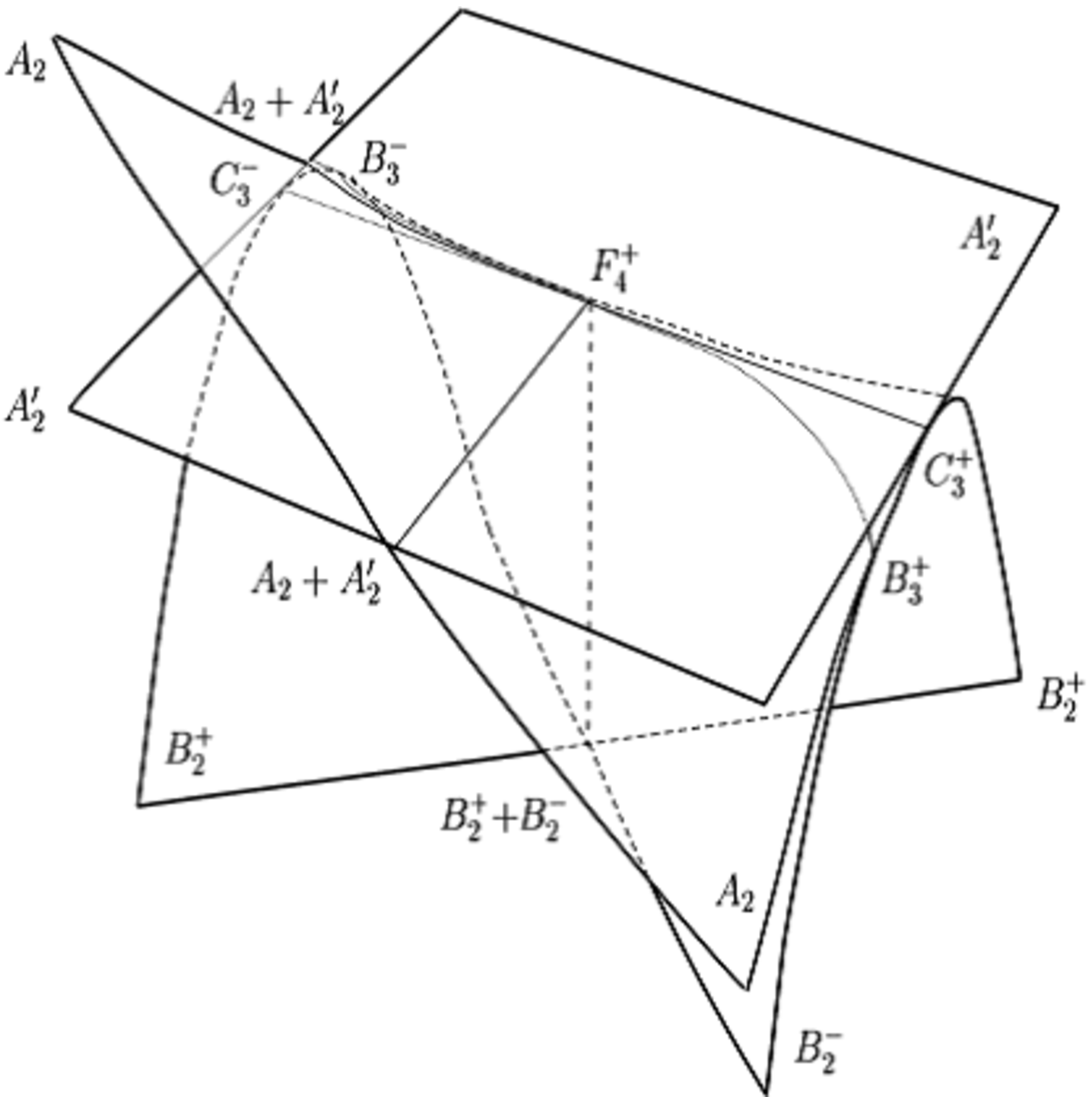}
&
\includegraphics*[width=7cm,height=7cm]{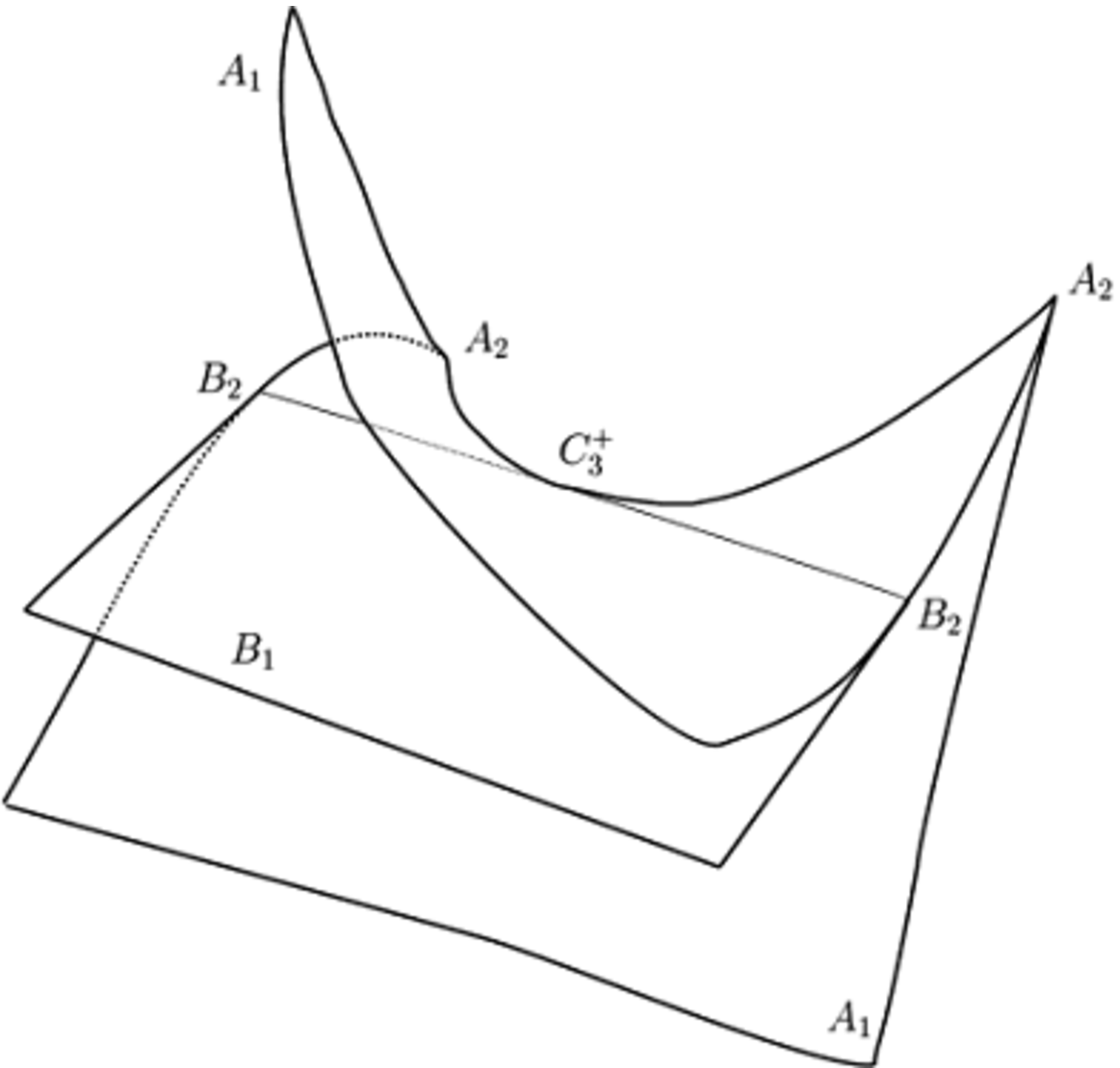}
\end{array}\]
\caption{the caustic $F^+_4$ and the wavefront $C^+_3$}
\end{figure}%

This paper consists of three parts.
In Part I, we recall stabilities of unfoldings under the equivalence relations {\it
reticular ${\cal P}$-${\cal R}^+$-equivalence} and 
{\it reticular ${\cal P}$-${\cal K}$-equivalence}.
They work as the equivalence relations of generating families of reticular Lagrangian, Legendrian maps respectively. 
In part \ref{Lagmap:part} and \ref{Legmap:part} we recall that the equivalence relations of 
reticular Lagrangian, Legendrian maps 
are reduced to that of their generating families.
In part \ref{Laggen:part} we shall study that the genericities of reticular Lagrangian, Legendrian maps 
are reduced to that of their generating families.
\part{Stability of unfoldings}

\section{Preliminaries}
\quad 
Let ${\mathbb H}^r=\{ (x_1,\ldots,x_r)\in {\mathbb R}^r|x_1\geq 0,\ldots,x_r\geq 0\}$ 
be an $r$-corner.
We denote by ${\cal E}(r;k_1,r;$ $k_2)$ the set of all germs at $0$ in 
${\mathbb H}^r\times {\mathbb R}^{k_1}$ of
smooth maps ${\mathbb H}^r\times {\mathbb R}^{k_1} \rightarrow 
{\mathbb H}^r\times {\mathbb R}^{k_2}$ and set ${\mathfrak M}(r;k_1,r;k_2)=
\{ f\in {\cal
E}(r;k_1,r;k_2)|f(0)=0 \}$.
We denote ${\cal E}(r;k_1,k_2)$ for ${\cal E}(r;k_1,0;k_2)$ and 
denote ${\mathfrak M}(r;k_1,k_2)$ for ${\mathfrak M}(r;k_1,0;k_2)$.

 If $k_2=1$ we write simply ${\cal E}(r;k)$ for 
${\cal E}(r;k,1)$
and ${\mathfrak M}(r;k)$ for ${\mathfrak M}(r;k,1)$. 
We also write ${\cal E}(k)$ for ${\cal E}(0;k)$ and ${\mathfrak M}(k)$
for ${\mathfrak M}(0;k)$.
Then ${\cal E}(r;k)$ is an ${\mathbb R}$-algebra in the usual
way and ${\mathfrak M}(r;k)$ is its unique maximal ideal. 

Let
$I_r=\{ 1,2,\ldots,r \}$ and $(x,y)=(x_1,\ldots,x_r,y_1,\ldots,y_k)$ be a fixed
coordinate system of $({\mathbb H}^r\times {\mathbb R}^k,0)$. 
We denote by 
${\cal B}(r;k)$ the group of
diffeomorphism germs on $({\mathbb H}^r\times {\mathbb R}^{k},0)$ of the form:
\[ \phi(x,y)=(x_1\phi_1^1(x,y),\ldots,x_r\phi_1^r(x,y),\phi_2^1(x,y),\ldots,\phi_2^k(x,y)
). \]
We also denote by 
${\cal B}_n(r;k+n)$ the group of
diffeomorphism germs on 
$({\mathbb H}^r\times {\mathbb R}^{k+n},0)$ of the form:
\[ \phi(x,y,u)=(x_1\phi_1^1(x,y,u),\ldots,x_r\phi_1^r(x,y,u),\phi_2^1(x,y,u),\ldots,\phi_2^k(x,y,u)
,\phi_3^1(u),\ldots,\phi_3^n(u)). \]


We denote $J^l(r+k,1)$ the set of $l$-jets at $0$ of germs in ${\mathfrak M}(r;k)$ 
and let $\pi_l:{\mathfrak M}(r;k)\rightarrow J^l(r+k,1)$ be the natural projection. 
We denote $j^lf(0)$ the $l$-jet of $f\in {\mathfrak M}(r;k)$. 
We also denote $\phi(x,y,u)=(x\phi_1(x,y,u),\phi_2(x,y,u),\phi_3(u))$ and denote other notations analogously.
In this paper all maps and all map germs are supposed to be smoothly. 
\section{Reticular ${\cal P}$-${\cal R}{}^+$-stability of unfoldings}
\quad 
We recall the stabilities of unfoldings under  
{\it the reticular ${\cal R}{}^+$-equivalence} which is developed in 
\cite{retLag}.
In order to distinguish equivalence relations between function germs in 
${\mathfrak M}(r;k)$
and their unfoldings, we denote this equivalence relation by
{\it the reticular ${\cal P}$-${\cal R}{}^+$-equivalence} in this paper.

We say that $f,g\in{\cal E}(r;k)$ are {\it reticular ${\cal R}$-equivalent} if
there exists $\phi\in{\cal B}(r;k)$ 
such that $g=f\circ \phi$. 

We say that function germs $f(x,y_1,\ldots,y_{k_1})
\in {\mathfrak M}(r;k_1)^2$ 
and $g(x,y_1,\ldots,y_{k_2})\in {\mathfrak M}(r;k_2)^2$ are {\it stably 
 reticular ${\cal R}$-equivalent} if $f$ and $g$ are reticular 
${\cal R}$-equivalent 
after additions of non-degenerate quadratic forms in the variables $y$. 
%
%

\vspace{2mm}
We say that a function germ $f\in {\mathfrak M}(r;k)$ is {\it reticular 
${\cal R}$-$l$-determined} 
if all function germ which has the same $l$-jet of $f$ is 
reticular ${\cal R}$-equivalent to $f$. 
\begin{lem}{\rm (cf., \cite[Lemma 4.2]{retLag})}\label{findet:lm}
Let $f(x,y)\in {\mathfrak M}(r;k)$ and let 
\[ {\mathfrak M}(r;k)^{l+1}\subset {\mathfrak M}(r;k)(\langle x_1 \frac{\partial f}{\partial
x_1},\ldots,x_r \frac{\partial f}{\partial x_r}\rangle +{\mathfrak M}(r;k)\langle
\frac{\partial f}{\partial y_1},\ldots,\frac{\partial f}{\partial y_k} \rangle )
+{\mathfrak M}(r;k)^{l+2},\]
then $f$ is reticular ${\cal R}$-$l$-determined.
Conversely if $f(x,y)\in {\mathfrak M}(r;k)$ be reticular ${\cal R}$-$l$-determined, then 
\[ {\mathfrak M}(r;k)^{l+1}\subset \langle x_1 \frac{\partial f}{\partial x_1},\ldots,x_r
\frac{\partial f}{\partial x_r}\rangle_{ {\cal E}(r;k) }  +{\mathfrak M}(r;k)\langle
\frac{\partial f}{\partial y_1},\ldots,\frac{\partial f}{\partial y_k}
\rangle. \]
\end{lem}   

We denote $x\frac{\partial f}{\partial x}$ for $(x_1 \frac{\partial f}{\partial x_1},\ldots,x_r)
\frac{\partial f}{\partial x_r}$ and $\frac{\partial f}{\partial y}$ for 
$(\frac{\partial f}{\partial y_1},\ldots,\frac{\partial f}{\partial y_k})$, and 
denote other notations analogously.

\vspace{5mm}
\quad 
Let $F(x,y,u)\in {\mathfrak M}(r;k+n_1),\ G(x,u,v)\in {\mathfrak M}(r;k+n_2)$ 
be unfoldings of $f\in {\mathfrak M}(r;k)$. 
We say that $G$ {\it is reticular ${\cal P}$-${\cal R}^+$-$f$-induced from} $F$  
if there exist 
$\Phi\in{\mathfrak M}(r;k+n_2,r;k+n_1)$ and
$\alpha \in {\mathfrak M}(n_2)$ 
satisfying the following conditions: \\
(1) $\Phi(x,y,0)=(x,y,0)$ for all $(x,y)\in ({\mathbb H}^r\times {\mathbb R}^k,0)$,\\
(2) $\Phi$ can be written in the form: 
$\Phi(x,y,v)=(x\phi_1(x,y,v),\phi_2(x,y,v),\phi_3(v))$,\\
(3) $G(x,y,v)=F\circ\Phi(x,y,v)+\alpha(v)$ for all 
    $(x,y,v)\in ({\mathbb H}^r\times {\mathbb R}^{k+n_2},0)$. \vspace{2mm}

We say that $F,G\in{\cal E}(r;k+n)$ are {\it reticular ${\cal P}$-${\cal R}^{(+)}$-equivalent} if
there exist $\Phi\in{\cal B}_n(r;k+n)$ 
(and $\alpha \in {\cal E}(n)$) such that $G=F\circ \Phi(+\alpha)$. 
We call $(\Phi,\alpha)$ a reticular ${\cal P}$-${\cal R}^{(+)}$-isomorphism 
from $F$ to $G$.
%
%
\begin{dfn}{\rm 
Here we recall the definitions of  several stabilities of unfoldings under the 
reticular ${\cal P}$-${\cal R}^+$-equivalence.  
Let $F(x,y,u)\in {\mathfrak M}(r;k+n)$ be an unfolding of 
$f(x,y)\in {\mathfrak M}(r;k)$.

\vspace{3mm}
%
We say that $F$ is {\it reticular 
${\cal P}$-${\cal R}^+$-stable} if the following condition holds: For any neighbourhood $U$ of
$0$ in ${\mathbb R}^{r+k+n}$ and any representative $\tilde{F} \in C^\infty 
(U,{\mathbb R})$ of $F$, 
there exists a neighbourhood $N_{\tilde{F}}$ of $\tilde{F}$ in $C^\infty$-topology
 such that for any element $\tilde{G} \in N_{\tilde{F}}$
the germ $\tilde{G}|_{{\mathbb H}^r\times {\mathbb R}^{k+n}}$ at 
$(0,y_0,u_0)$ is reticular ${\cal P}$-${\cal R}{}^+$-equivalent to $F$
for some $(0,y_0,u_0)\in U$.

\vspace{3mm}
%
We say that $F$ is {\it reticular
${\cal P}$-${\cal R}{}^+$-versal} if all unfolding of $f$ is  reticular ${\cal P}$-${\cal R}{}^+$-$f$-induced from $F$.
\vspace{3mm}

%
We say that $F$ is {\it reticular ${\cal P}$-${\cal R}{}^+$-infinitesimally versal} if 
\[ 
{\cal E}(r;k)=\langle x \frac{\partial f}{\partial x},
\frac{\partial f}{\partial
y}\rangle_{ {\cal E}(r;k) }+
\langle 1,\frac{\partial F}{\partial u}|_{u=0}\rangle_{{\mathbb R}}.
\] 
\vspace{3mm}

%
\quad 
We say that $F$ is {\it reticular ${\cal P}$-${\cal R}^+$-infinitesimally stable} if
\[ {\cal E}(r;k+n)  = \langle x \frac{\partial
 F}{\partial x},\frac{\partial
 F}{\partial y}\rangle_{ {\cal
 E}(r;k+n) }+\langle 1,\frac{\partial F}{\partial u}\rangle_{{\cal E}(n)}. 
\]
\vspace{3mm}

%
\quad 
We say that $F$ is {\it reticular ${\cal P}$-${\cal R}^+$-homotopically stable} if for any
smooth path-germ $({\mathbb R},0)\rightarrow {\cal E}(r;k+n),t\mapsto F_t$ with
$F_0=F$, there exists a smooth path-germ $({\mathbb R},0)\rightarrow {\cal
B}_n(r;k+n)\times {\cal E}(n),t\mapsto (\Phi_t,\alpha_t)$ with
$(\Phi_0,\alpha_0)=(id,0)$ such that each $(\Phi_t,\alpha_t)$ is a reticular
${\cal P}$-${\cal R}{}^+$-isomorphism from $F$ to $F_t$, that is
$F_t=F\circ \Phi_t+\alpha_t$ for $t$ around $0$.
}\end{dfn}
%
%
\begin{tth}{\rm (cf., \cite[Theorem 4.5]{retLag})}
Let $F\in {\mathfrak M}(r;k+n)$ be an unfolding of $f\in {\mathfrak M}(r;k)$.
Then the following are equivalent. \\
{\rm (1)} $F$ is reticular ${\cal P}$-${\cal R}{}^+$-stable.\\
{\rm (2)} $F$ is reticular ${\cal P}$-${\cal R}{}^+$-versal.\\
{\rm (3)} $F$ is reticular ${\cal P}$-${\cal R}{}^+$-infinitesimally versal. \\
{\rm (4)} $F$ is reticular ${\cal P}$-${\cal R}{}^+$-infinitesimally stable. \\
{\rm (5)} $F$ is reticular ${\cal P}$-${\cal R}{}^+$-homotopically stable. 
\end{tth}

For a function germ $f(x,y) \in {\mathfrak M}(r;k)$, 
if $1,a_1,\ldots,a_n\in {\cal E}(r;k)$ is a representative of a basis of the 
vector space 
\[ {\cal E}(r;k)/ \langle x \frac{\partial f}{\partial x},
\frac{\partial f}{\partial y}\rangle_{ {\cal E}(r;k)}, \]
then 
the function germ $f+a_1u_1+\cdots +a_nu_n\in {\mathfrak M}(r;k+n)$ is a 
reticular ${\cal P}$-${\cal R}{}^+$-stable unfolding of $f$.
We call $n$ the reticular ${\cal R}^+$-codimension of $f$.

\vspace{3mm}
We call a function germ $f\in {\mathfrak M}(r;k)$ is {\it ${\cal R}$-simple} if 
the following holds:
For a sufficiently higher integer $l$,
there exists a neighbourhood $N$ of $j^lf(0)$ in $J^l(r+k,1)$ such that 
$N$ intersects finite ${\cal R}$-orbits. 
By \S 17.4\cite[p.279]{arnold:text} we have that:
\begin{tth}\label{simpleR:th}
 An ${\cal R}$-simple function germ in ${\mathfrak M}(1;k)^2$ is 
stably ${\cal R}$-equivalent to one of the following function germ:
\[
B^{\pm}_l:\pm x^l\ (l\geq 2),\ \ \ C^{\pm}_l:xy\pm y^{l}\ (l\geq 3),\ \ \ 
F^{\pm}_4: \pm x^2+y^3.\]
 \end{tth} 
\section{Reticular ${\cal P}$-${\cal K}$-stability of unfoldings}
\quad 
We recall the stabilities of  unfoldings under 
{\it the reticular ${\cal K}$-equivalence} which is developed in 
\cite{retLeg}.
In this paper we denote this equivalence relation by 
{\it the reticular ${\cal P}$-${\cal K}$-equivalence}.

We say that $f,g\in{\cal E}(r;k)$ are {\it reticular ${\cal K}$-equivalent} if
there exist $\phi\in{\cal B}(r;k)$ and a unit $a\in {\cal E}(r;k)$
such that $g=a\cdot f\circ \phi$. 

We say that function germs $f(x,y_1,\ldots,y_{k_1})
\in {\mathfrak M}(r;k_1)^2$ 
and $g(x,y_1,\ldots,y_{k_2})\in {\mathfrak M}(r;k_2)^2$ are {\it stably 
 reticular ${\cal K}$-equivalent} if $f$ and $g$ are reticular 
${\cal K}$-equivalent 
after additions of non-degenerate quadratic forms in the variables $y$. 

\vspace{3mm}
\quad
We say that a function germ $f\in {\mathfrak M}(r;k)$ is {\it reticular ${\cal K}$-$l$-determined} 
if all function germ which has the same $l$-jet of $f$ is 
reticular ${\cal K}$-equivalent to $f$. 
\begin{lem}{\rm (cf., \cite[Lemma 6.2]{retLeg})}\label{findetc:lm}
Let $f(x,y)\in {\mathfrak M}(r;k)$ and let 
\[ {\mathfrak M}(r;k)^{l+1}\subset {\mathfrak M}(r;k)(\langle f,
x \frac{\partial f}{\partial x}\rangle +{\mathfrak M}(r;k)\langle
\frac{\partial f}{\partial y}\rangle )
+{\mathfrak M}(r;k)^{l+2},\]
then $f$ is reticular ${\cal K}$-$l$-determined. 
Conversely if $f(x,y)\in {\mathfrak M}(r;k)$ be reticular ${\cal K}$-$l$- determined, then 
\[ {\mathfrak M}(r;k)^{l+1}\subset \langle f,
x \frac{\partial f}{\partial x}\rangle_{ {\cal E}(r;k) }  
+{\mathfrak M}(r;k)\langle \frac{\partial f}{\partial y}\rangle. \]
\end{lem}   
\quad 
Let $F(x,y,u)\in {\mathfrak M}(r;k+n_1),\ G(x,y,v)\in {\mathfrak M}(r;k+n_2)$ 
be unfoldings of $f(x,y)\in {\mathfrak M}(r;k)$. 
We say that $G$ {\it is reticular ${\cal P}$-${\cal K}$-$f$-induced from} $F$  
if there exist 
$\Phi\in{\mathfrak M}(r;k+n_2,r;k+n_1)$ and
$\alpha \in {\cal E}(r;k+n_2)$
satisfying the following conditions: \\
(1) $\Phi(x,y,0)=(x,y,0)$ and $\alpha(x,y,0)=1$ for all $(x,y)\in ({\mathbb H}^r\times {\mathbb R}^k,0)$,\\
(2) $\Phi$ can be written in the form:
$\Phi(x,y,v)=(x\phi_1(x,y,v),
\phi_2(x,y,v),\phi_3(v)),$\\
(3) $G(x,y,v)=\alpha(x,y,v)\cdot F\circ\Phi(x,y,v)$ for all 
    $(x,y,v)\in ({\mathbb H}^r\times {\mathbb R}^{k+n_2},0)$
\vspace{2mm} 

We say that $F,G\in{\cal E}(r;k+n)$ are {\it reticular ${\cal P}$-${\cal K}$-equivalent} if
there exist $\Phi\in{\cal B}_n(r;k+n)$ and a unit $\alpha\in {\cal E}(r;k+n)$
such that $G=\alpha\cdot F\circ \Phi$.
We call $(\Phi,\alpha)$ a reticular ${\cal P}$-${\cal K}$-isomorphism from $F$ to $G$. 
%
%
\begin{dfn}{\rm 
Here we recall the definitions of several stabilities of unfoldings under the 
reticular ${\cal P}$-${\cal K}$-equivalence.  
Let $F(x,y,u)\in {\mathfrak M}(r;k+n)$ be an unfolding of 
$f(x,y)\in {\mathfrak M}(r;k)$.

\vspace{3mm}
%
We say that $F$ is {\it reticular 
${\cal P}$-${\cal K}$-stable} if the following condition holds: For any neighbourhood $U$ of
$0$ in ${\mathbb R}^{r+k+n}$ and any representative $\tilde{F} \in C^\infty 
(U,{\mathbb R})$ of $F$, there exists a neighbourhood $N_{\tilde{F}}$ of $\tilde{F}$ 
in $C^\infty$-topology such that for any element $\tilde{G} \in N_{\tilde{F}}$
the germ $\tilde{G}|_{{\mathbb H}^r\times {\mathbb R}^{k+n}}$ at 
$(0,y_0,u_0)$ is reticular ${\cal P}$-${\cal K}$-equivalent to $F$
for some $(0,y_0,u_0)\in U$.

\vspace{3mm}
%
We say that $F$ is {\it reticular
${\cal P}$-${\cal K}$-versal} if all unfolding of $f$ is reticular ${\cal P}$-${\cal K}$-$f$-induced 
from $F$.

\vspace{3mm}
%
We say that $F$ is {\it reticular ${\cal P}$-${\cal K}$-infinitesimally versal} if 
\[ 
{\cal E}(r;k)=\langle f,x \frac{\partial f}{\partial x},
\frac{\partial f}{\partial y}\rangle_{ {\cal E}(r;k) }+
\langle \frac{\partial F}{\partial
u}|_{u=0}\rangle_{{\mathbb R}}.
\] 

\vspace{3mm}
%
We say that $F$ is {\it reticular ${\cal P}$-${\cal K}$-infinitesimally stable} if
\[
{\cal E}(r;k+n)= \langle F,x \frac{\partial
 F}{\partial x},\frac{\partial
 F}{\partial y}\rangle_{ {\cal
 E}(r;k+n) }+\langle \frac{\partial F}{\partial u}\rangle_{{\cal E}(n)}. 
\]

\vspace{3mm}
%
\quad 
We say that $F$ is {\it reticular ${\cal P}$-${\cal K}$-homotopically stable} if for any
smooth path-germ $({\mathbb R},0)\rightarrow {\cal E}(r;k+n),t\mapsto F_t$ with
$F_0=F$, there exists a smooth path-germ $({\mathbb R},0)\rightarrow {\cal
B}_n(r;k+n)\times {\cal E}(r;k+n),t\mapsto (\Phi_t,\alpha_t)$ with
$(\Phi_0,\alpha_0)=(id,1)$ such that each $(\Phi_t,\alpha_t)$ is a reticular
${\cal P}$-${\cal K}$-isomorphism from $F$ to $F_t$, that is 
$F_t=\alpha_t\cdot F\circ \Phi_t$ for $t$ around $0$.
}\end{dfn}
%
%
\begin{tth}{\rm (cf., \cite[Theorem 6.5]{retLeg})}
Let $F\in {\mathfrak M}(r;k+n)$ be an unfolding of $f\in {\mathfrak M}(r;k)$. 
Then the following are equivalent. \\
{\rm (1)} $F$ is reticular ${\cal P}$-${\cal K}$-stable.\\
{\rm (2)} $F$ is reticular ${\cal P}$-${\cal K}$-versal.\\
{\rm (3)} $F$ is reticular ${\cal P}$-${\cal K}$-infinitesimally versal. \\
{\rm (4)} $F$ is reticular ${\cal P}$-${\cal K}$-infinitesimally stable. \\
{\rm (5)} $F$ is reticular ${\cal P}$-${\cal K}$-homotopically stable. 
\end{tth}

For a function germ $f(x,y) \in {\mathfrak M}(r;k)$, 
if $a_1,\ldots,a_n\in {\cal E}(r;k)$ is a representative of a basis of the 
vector space 
\[ {\cal E}(r;k)/ \langle f,x \frac{\partial f}{\partial x},\frac{\partial f}{\partial
y}\rangle_{ {\cal E}(r;k)}, \]
then 
the function germ $f+a_1u_1+\cdots +a_nu_n\in {\mathfrak M}(r;k+n)$ is a 
reticular ${\cal P}$-${\cal K}$-stable unfolding of $f$.
We call $n$ the reticular ${\cal K}$-codimension of $f$.

\vspace{3mm}
We call a function germ $f\in {\mathfrak M}(r;k)$ is 
{\it ${\cal K}$-simple} if 
the following holds:
For a sufficiently higher integer $l$,
there exists a neighbourhood $N$ of $j^lf(0)$ in $J^l(r+k,1)$ such that 
$N$ intersects finite ${\cal K}$-orbits. 
\begin{tth}\label{simpleK:th}
 A ${\cal K}$-simple function germ in ${\mathfrak M}(1;k)^2$ is 
stably ${\cal K}$-equivalent to one of the following function germ:
\[
B_l:x^l\  (l\geq 2),\ \ \ C^{\vare}_l:xy+\vare y^{l}\ (\vare^{l-1}=1,
l\geq 3),\ \ \ 
F_4: x^2+y^3.\]
 \end{tth} 
\part{Reticular Lagrangian maps}\label{Lagmap:part}
\section{Symplectic regular $r$-cubic configurations}
\quad 
Let $(q,p)=(q_1,\ldots,q_n,p_1,\ldots,p_n)$ be a canonical coordinate system of $(T^* {\mathbb R}^n,0)$ and 
$\pi:(T^* {\mathbb R}^n,0)\rightarrow ({\mathbb R}^n,0)$ be the cotangent 
bundle equipped with the canonical symplectic structure $dp\wedge dq$.  
We define 
\[ L_\sigma^0=
\{(q,p)\in (T^* {\mathbb R}^n,0)|q_\sigma=p_{I_r-\sigma}=q_{r+1}=\cdots=q_n=0, 
q_{I_r-\sigma}\geq 0 \} \]
for $\sigma\subset I_r$.

\begin{dfn}{\rm
Let $\{L_\sigma \}_{\sigma \subset I_r}$ be a family of $2^r$
 Lagrangian submanifold germs on $(T^* {\mathbb R}^n,0)$. 
Then $\{L_\sigma \}_{\sigma \subset I_r}$ is called 
{\it a symplectic regular $r$-cubic configuration} if there exists a 
symplectic diffeomorphism germ $S$ on $(T^* {\mathbb R}^n,0)$
 such that 
$L_\sigma=S(L^0_\sigma)$ for all $\sigma \subset I_r$.
}\end{dfn}
\section{Reticular Lagrangian maps and their generating $ $ families}
\quad
We introduce a main result in \cite{retLag} about 
the relations of reticular Lagrangian maps and their generating families.\\
 
Let ${\mathbb L}=\{(q,p)\in T^* {\mathbb R}^n |q_1p_1=\cdots 
=q_rp_r=q_{r+1}=\cdots=q_n=0,q_{I_r}\geq 0  \}$ be the representative 
as a germ of the union of $L^0_\sigma$ for $\sigma\subset I_r$.
We call a map germ 
\[ ({\mathbb L},0) \stackrel{i}{\longrightarrow} (T^* {\mathbb R}^n,0) \stackrel{\pi}{\longrightarrow} ({\mathbb R}^n,0) \]
{\it a reticular Lagrangian map} if there exists a symplectic diffeomorphism 
germ
$S$ on $(T^* {\mathbb R}^n,0)$ such that $i=S|_{{\mathbb L}}$. 
We call $i$ {\it the reticular Lagrangian embedding} of $\pi\circ i$. 
We call $S$ {\it an extension} of $i$ and call 
$\{ i(L^0_\sigma) \}_{ \sigma \subset I_r }$  {\it the symplectic regular
$r$-cubic configuration} associated with $\pi\circ i$.

{\bf Caustics}: 
Let $\pi\circ i$ be a reticular Lagrangian map. 
Let $C_\sigma$ be the caustic  of the Lagrangian map  
$\pi\circ i|_{L^0_\sigma}$ for $\sigma\subset I_r$,
that is, the set of  critical values of $\pi\circ i|_{L^0_\sigma}$.
Let  
$Q_{\sigma,\tau}=\pi\circ i(L^0_\sigma \cap L^0_\tau)$ 
for $\sigma\neq\tau\subset I_r$. 
We define the {\it caustic} of $\pi\circ i$ by 
\[ \bigcup_{\sigma \subset I_r} C_\sigma \cup \bigcup_{\sigma\neq \tau}
Q_{\sigma,\tau}. \]
We remark that for $\tau_1,\tau_2\subset I_r\ (\tau_1\neq\tau_2)$
we have $Q_{\tau_1,\tau_2}\subset Q_{\sigma,\sigma\cup \{ i\}}$, where 
$\sigma=\tau_1\cap\tau_2$ and $i$ be any element of 
$(\tau_1-\sigma)\cup (\tau_2-\sigma)$. 
This means that $\bigcup_{\sigma\neq \tau} Q_{\sigma,\tau}$ 
is equal to the union of $Q_{\sigma,\tau}$ for 
$\sigma\subset \tau\subset I_r,\ \#(\tau-\sigma)=1$. 
For example, in the case $r=2$ we have 
\[ \bigcup_{\sigma\neq \tau} Q_{\sigma,\tau}\ = \ 
Q_{\emptyset,1}\cup Q_{\emptyset,2}\cup Q_{1,\{ 1,2\}}\cup Q_{2,\{1,2\}}. \]

\vspace{3mm}
A function germ $F(x,y,q)\in {\mathfrak M}^2(r;k+n)$ is called 
{\it $S$-non-degenerate } if 
\[ x_1,\ldots,x_r,\frac{\partial F}{\partial x_1},\ldots,\frac{\partial F}{\partial x_r},\frac{\partial F}{\partial y_1},\ldots,\frac{\partial F}{\partial y_k} \]
are independent on $({\mathbb H}^k\times {\mathbb R}^{k+n},0)$, that is
\[ \mbox{rank}\left( 
\begin{array}{cc}
\displaystyle{\frac{\partial^2 F}{\partial x\partial y}} & 
\displaystyle{\frac{\partial^2 F}{\partial x\partial u}} \\
\displaystyle{\frac{\partial^2 F}{\partial y\partial y}} & 
\displaystyle{\frac{\partial^2 F}{\partial y\partial u}} \\
\end{array}
\right)_0
=r+k.\]
\begin{dfn}{\rm 

Let $\{L_\sigma \}_{\sigma \subset I_r}$ be a 
symplectic regular $r$-cubic configuration in $(T^* {\mathbb R}^n,0)$ and 
$F(x,y,q)\in {\mathfrak M}(r;k+n)^2$ be a  function germ. 
We call $F$ {\it a generating family} of 
$\{L_\sigma \}_{\sigma \subset I_r}$ if the following conditions hold:\\
(1) $F$ is $S$-non-degenerate,\\
(2) 
$F|_{x_\sigma=0}$ is a generating family of $L_\sigma$ for 
$\sigma\subset I_r$, that is 
\[ L_\sigma=\{ (q,\frac{\partial F}{\partial q}(x,y,q))
  \in (T^*{\mathbb R}^n ,0)|
     x_\sigma=\frac{\partial F}{\partial x_{I_r-\sigma}}=
              \frac{\partial F}{\partial y}=0,x_{I_r-\sigma}\geq 0 \}. \]
We also call $F$ a generating family of a reticular Lagrangian map 
$\pi\circ i$
if $F$ is a generating family of 
the symplectic regular $r$-cubic configuration 
$\{i(L^0_\sigma)\}_{\sigma\subset I_r}$.
}\end{dfn}

We call 
a symplectic diffeomorphism 
germ 
$\phi$ on $(T^* {\mathbb R}^n,0)$ 
 {\it a reticular diffeomorphism} if  
$\phi(L^0_\sigma)=L^0_\sigma$ for $\sigma \subset I_r$. 
We say that reticular Lagrangian maps $\pi \circ i_1,\pi \circ i_2:(
{\mathbb L},0) \rightarrow (T^* {\mathbb R}^n,0) \rightarrow({\mathbb R}^n,0)$ are 
{\it Lagrangian equivalent} if there exist a reticular diffeomorphism $\phi$ and a Lagrangian equivalence $\Theta$ of $\pi$ such that the following diagram is commutative:
\[
\begin{array}{ccccc}
({\mathbb L},0) & \stackrel{i_1}{\longrightarrow} & (T^* {\mathbb R}^n,0) & \stackrel{\pi}{\longrightarrow} & ({\mathbb R}^n,0) \\
\phi|_{\mathbb L}  \downarrow \ \ \ \ & & \Theta \downarrow \ \ & & g \downarrow \ \ \\
({\mathbb L},0) & \stackrel{i_2}{\longrightarrow} & (T^* {\mathbb R}^n,0) & \stackrel{\pi}{\longrightarrow} & ({\mathbb R}^n,0), \\
\end{array}
\]
where $g$ is the diffeomorphism of the base space 
of $\pi$ induced from $\Theta$.

We remark that there is not the condition that a reticular diffeomorphism 
is a symplectic diffeomorphism in the definition in \cite{retLag}.
But a reticular diffeomorphism defined in \cite{retLag} consists of 
a restriction to ${\mathbb L}$ of 
compositions of two symplectic diffeomorphism and a Lagrangian equivalence,
it follows that a reticular diffeomorphism  in \cite{retLag} is automatically a 
 restriction to ${\mathbb L}$ of 
symplectic diffeomorphism.
So our definition is equivalent to that in \cite{retLag}.\\

We say that function germs $F(x,y_1,\ldots,y_{k_1},q)\in {\mathfrak M}(r;k_1+n)$ 
and $G(x,y_1,\ldots,y_{k_2},q)\in {\mathfrak M}(r;k_2+n)$ are {\it stably 
 reticular ${\cal P}$-${\cal R}^+$-equivalent} if $F$ and $G$ are reticular 
${\cal P}$-${\cal R}^+$-equivalent 
after additions of non-degenerate quadratic forms in the variables $y$. 
\begin{tth}\label{Laggene:th}{\rm (cf., \cite[Theorem 3.2]{retLag})}
{\rm (1)} For any reticular Lagrangian map $\pi \circ i$, there exists a function germ $F\in {\mathfrak M}(r;k+n)^2$ which is a generating family of $\pi\circ i$. \\
{\rm (2)} For any S-non-degenerate function germ $F\in {\mathfrak M}(r;k+n)^2$, there exists a reticular Lagrangian map of which $F$ is a generating family. \\
{\rm (3)} Two reticular Lagrangian maps are Lagrangian equivalent if and only if their generating families are stably reticular ${\cal P}$-${\cal R}^{+}$-equivalent. 
\end{tth}
\section{Stabilities of reticular Lagrangian  maps}
\quad
We introduce several stabilities of reticular Lagrangian maps:\\
{\bf Stability}:
For any open set $U$ in $T^*{\mathbb R}^n$, 
we denote $S(U,T^*{\mathbb R}^n)$ the space of symplectic embeddings from 
$U$ to $T^*{\mathbb R}^n$ with $C^\infty$-topology. 
We say that 
a reticular Lagrangian map $\pi \circ i:({\mathbb L}^0,0) \rightarrow (T^* {\mathbb
  R}^n,0) \rightarrow({\mathbb R}^n,0)$ is {\it stable} if the following
holds: For any extension $S$ of $i$ and any representative
$\tilde{S}\in S(U,T^*{\mathbb R}^n)$ of $S$, there exists a neighbourhood
$N_{\tilde{S}}$ of $\tilde{S}$ such that 
for any $\tilde{S'}\in N_{\tilde{S}}$ 
the reticular Lagrangian maps $\pi \circ \tilde{S}'_{w_0}|_{{\mathbb L}^0}$
and $\pi \circ i$ are Lagrangian equivalent for some
$w_0=(0,\ldots,0,p^0_{r+1},\ldots,p^0_n)\in U$,
where the symplectic diffeomorphism germ $\tilde{S}'_{w_0}$ on 
$(T^* {\mathbb R}^n,0)$ is defined by $w\mapsto \tilde{S'}(w+w_0)-\tilde{S'}(w_0)$. 
\vspace{2mm}\\
{\bf Homotopically stability}:
Let $\pi \circ i:({\mathbb L},0) \rightarrow (T^* {\mathbb R}^n,0)
\rightarrow({\mathbb R}^n,0)$ be a reticular Lagrangian map. 
A one-parameter family
of symplectic diffeomorphisms $\bar{S}:(T^* {\mathbb R}^n\times {\mathbb R},(0,0))\rightarrow
(T^* {\mathbb R}^n,0)((Q,P,t)\mapsto S_t(Q,P))$ 
 is called a {\it reticular Lagrangian
  deformation} of $i$ 
if $S_0|_{\mathbb L}=i_0$.
Let $\phi$ be a reticular
diffeomorphism on $(T^* {\mathbb R}^n,0)$. 
A map germ $\bar{\phi}:(T^* {\mathbb R}^n \times {\mathbb R},(0,0)) \rightarrow 
(T^* {\mathbb R}^n,0)((Q,P,t)\mapsto \phi_t(Q,P))$ is called a {\it one-parameter
  deformation of reticular diffeomorphisms} of $\phi$ if $\phi_0=\phi$ and 
$\phi_t$ 
 is a reticular diffeomorphism for $t$ around $0$. 
We say that a reticular Lagrangian map $\pi \circ i:({\mathbb
  L},0) \rightarrow (T^* {\mathbb R}^n,0) \rightarrow({\mathbb R}^n,0)$ is {\it
  homotopically stable} if for any reticular Lagrangian deformation 
$\bar{S}=\{ S_t \}$ of $i$
, there exists a one-parameter deformation of
reticular diffeomorphisms $\bar{\phi}=\{ \phi_t \}$ of $
id_{T^* {\mathbb R}^n}$ and
 a one-parameter family of Lagrangian equivalences
$\bar{\Theta}=\{ \Theta_t \}$ of $\pi$ with $\Theta_0=id_{T^* {\mathbb R}^n}$
such that $S_t=\Theta_t\circ S_0\circ \phi_t$ for $t$ around
$0$.
\vspace{2mm}\\
{\bf Infinitesimally stability}:
Let $S$ be a symplectic diffeomorphism on $(T^* {\mathbb R}^n,0)$.
We say that a vector field $v$ on $(T^* {\mathbb R}^n,0)$ along $S$ 
is {\it an infinitesimal symplectic transformation of} $S$ if 
there exists a  reticular Lagrangian  deformation
 $\bar{S}=\{S_t\}$ such that $S_0=S$ and
$\frac{dS_t}{dt}|_{t=0}=v$.
We say that a vector field $\xi$ on $(T^* {\mathbb R}^n,0)$ is 
{\it an infinitesimally reticular diffeomorphism} if there exists  
a one-parameter deformation  of
reticular diffeomorphisms  $\bar{\phi}=\{ \phi_t \}$ of
$id_{T^* {\mathbb R}^n}$ 
such that  $\frac{d\phi_t}{dt}|_{t=0}=\xi$.
We say that a vector field $\eta$ on  $(T^* {\mathbb R}^n,0)$ is
{\it an infinitesimally Lagrangian equivalence} if 
there exists a one-parameter family of Lagrangian equivalences
$\bar{\Theta}=\{\Theta_t\}$
of $\pi$ 
such that $\Theta_0=id_{T^* {\mathbb R}^n}$ and 
$\frac{d\Theta_t}{dt}|_{t=0}=\eta$.
We say that a reticular
Lagrangian map $\pi \circ i:({\mathbb L},0) \rightarrow (T^* {\mathbb R}^n,0)
\rightarrow({\mathbb R}^n,0)$ is  
 {\it infinitesimally stable} if for any
extension $S$ of $i$ and any infinitesimally symplectic transformation $v$ of $S$,
there exists an infinitesimally reticular diffeomorphism $\xi$ and 
an infinitesimally Lagrange equivalence $\eta$ such that $v=S_*\xi+\eta\circ S$.

\vspace{3mm}
We say that a function germ $H$ on $(T^*{\mathbb R}^n,0)$ is {\it fiber preserving} if 
$H$ has the form $H(q,p)=\sum_{j=1}^nh_j(q)p_j+h_0(q)$.

\vspace{3mm}
We recall the following theorem which is proved in \cite{retLag}.
\begin{tth}\label{stable:thlag}{\rm (cf., \cite[p.587 Theorem 5.5]{retLag})}
Let $\pi \circ i:({\mathbb L},0) \rightarrow (T^* {\mathbb R}^n,0) \rightarrow({\mathbb
  R}^n,0)$ be a reticular Lagrangian map with a generating family
$F(x,y,q)\in {\mathfrak M}(r;k+n)^2$. Then the following
are equivalent.\\  
{\rm (u)} $F$ is a reticular ${\cal R}^+$-stable unfolding of $F|_{q=0}$. \\
{\rm (hs)} $\pi \circ i$ is homotopically stable. \\
{\rm (is)} $\pi \circ i$ is infinitesimally stable. \\
{\rm (a)} For any function germ $f$ on $(T^* {\mathbb R}^n,0)$, there exists a fiber preserving 
function germ $H$ on $ $ $(T^* {\mathbb R}^n,0)$ such that $f\circ i=H\circ i$. \\
{\rm (s)} $\pi \circ i$ is  stable. 
\end{tth}

The definition of the infinitesimally stability of reticular 
Lagrangian maps in  \cite{retLag} seem to be different from 
our one.
But these are equivalent by Lemma \ref{infsta:Laglem}.
The assertion that (u), (hs),(is) and (a) are all equivalent is proved in \cite{retLag}.
But the proof of (u)$\Leftrightarrow$(s) is slightly complicated.
So we shall prove this assertion by another way 
in  Section \ref{GeneLag}.
\part{Reticular Legendrian maps}\label{Legmap:part}
\section{Contact regular $r$-cubic configurations}
\quad 
Let $J^1({\mathbb R}^n,{\mathbb R})$ be the $1$-jet bundle of functions 
in $n$-variables which may be considered as ${\mathbb R}^{2n+1}$ with 
a natural coordinate system $(q,z,p)=(q_1,\ldots,q_n,z,p_1,\ldots,p_n)$,
where $q$ is a coordinate system of ${\mathbb R}^n$.
We equip the contact structure on $J^1({\mathbb R}^n,{\mathbb R})$
defined by the canonical $1$-form $\theta=dz-\sum_{j=1}^np_jdq_j$.
We have a natural projection 
$\tilde{\pi}:J^1({\mathbb R}^n,{\mathbb R})\rightarrow {\mathbb R}^n\times {\mathbb R}$
by $\tilde{\pi}(q,z,p)=(q,z)$.

\vspace{2mm}

Let $\tilde{L}^0_\sigma=
\{ (q,z,p) \in (J^1({\mathbb R}^n,{\mathbb R}),0)|q_\sigma=p_{I_r-\sigma}=
q_{r+1}=\cdots=q_n=z=0,q_{I_r-\sigma}\geq 0 \}$ for each $\sigma \subset I_r$.
\begin{dfn}\label{config}{\rm
Let $\{L_\sigma \}_{\sigma \subset I_r}$ be a family of $2^r$
 Legendrian submanifold germs on $(J^1({\mathbb R}^n,{\mathbb R}),$ $0)$. 
Then $\{L_\sigma \}_{\sigma \subset I_r}$ is called 
{\it a contact regular $r$-cubic configuration} if there exists a 
contact diffeomorphism germ $C$ on $(J^1({\mathbb R}^n,{\mathbb R}),0)$
 such that 
$L_\sigma=C(\tilde{L}^0_\sigma)$ for all $\sigma \subset I_r$.

}\end{dfn}
\section{Reticular Legendrian maps and their generating $ $ families}
\quad
We introduce a main result in \cite{retLeg} about 
the relations of reticular Legendrian maps and their generating families.\\

Let $\tilde{\mathbb L}=\{(q,z,p) \in J^1({\mathbb R}^n,{\mathbb R})|
q_1p_1=\cdots=q_rp_r=q_{r+1}=\cdots=q_n=z=0,q_{I_r}\geq 0 \}$ be a 
representative as a germ 
of the union of $\tilde{L}^0_\sigma$ for all $\sigma\subset I_r$.
We call a map germ 
\[ (\tilde{\mathbb L},0) \stackrel{i}{\longrightarrow} (J^1({\mathbb R}^n,
{\mathbb R}),0) \stackrel{\tilde{\pi}}{\longrightarrow} ({\mathbb R}^n\times 
{\mathbb R},0) \]
{\it a reticular Legendrian map} if there exists a contact diffeomorphism germ
$C$ on $(J^1({\mathbb R}^n,{\mathbb R}),0)$ such that $i=C|_{\tilde{\mathbb L}}$. 
We call $i$ {\it the reticular Legendrian embedding of} $\tilde{\pi} \circ i$.

{\bf Wavefronts}: Let $\tilde{\pi}\circ i$ be a reticular Legendrian map.
We define the {\it wavefront} of $\tilde{\pi}\circ i$  by 
the union of $\tilde{\pi}\circ i(L^0_\sigma)$ for all $\sigma \subset I_r$.

\vspace{2mm}

A function germ $F(x_1,\ldots,x_r,y_1,\ldots,y_k,q_1,\ldots,q_n,z)
\in {\mathfrak M}(r;k+n+1)$ is called {\it C-non-degenerate } 
if $\frac{\partial F}{\partial x}(0)=
\frac{\partial F}{\partial y}(0)=0$ and
\[ x_1,\ldots,x_r,F,\frac{\partial F}{\partial x_1},\ldots,
\frac{\partial F}{\partial x_r},\frac{\partial F}{\partial y_1},
\ldots,\frac{\partial F}{\partial y_k} \]
are independent on $({\mathbb H}^r\times {\mathbb R}^{k+n+1},0)$, that is
\[ {\rm rank}\left( 
\begin{array}{ccc}
\displaystyle{\frac{\partial F}{\partial y}} & 
\displaystyle{\frac{\partial F}{\partial q} } &
\displaystyle{\frac{\partial F}{\partial z}}\\
\displaystyle{\frac{\partial^2 F}{\partial x\partial y} } & 
\displaystyle{\frac{\partial^2 F}{\partial x\partial q} } &
\displaystyle{\frac{\partial^2 F}{\partial x\partial z} }\\
\displaystyle{\frac{\partial^2 F}{\partial y\partial y} } & 
\displaystyle{\frac{\partial^2 F}{\partial y\partial q} } &
\displaystyle{\frac{\partial^2 F}{\partial y\partial z} }\\
\end{array}
\right)_0
=r+k+1.\]
\begin{dfn}{\rm 
Let $\{ L_\sigma \}_{ \sigma \subset I_r }$ be a contact regular 
$r$-cubic configuration in $(J^1({\mathbb R}^n,{\mathbb R}),0)$. 
Then a function germ $F(x,y,q,z)\in {\mathfrak M}(r;k+n+1)$ is called 
{\it a generating 
family} of $\{ L_\sigma \}_{ \sigma \subset I_r }$ if 
the following conditions hold:\\
(1) $F$ is C-non-degenerate,\\
(2) For each $\sigma \subset I_r$, the function germ  $F|_{x_\sigma=0}$ 
is a generating 
family of $L_\sigma$, that is
\[ L_\sigma=\{ (q,z,
\frac{\partial F}{\partial q}/
(-\frac{\partial F}{\partial z}))\in (J^1({\mathbb R}^n,{\mathbb R}),0)|\ 
x_\sigma=\frac{\partial F}{\partial x_{I_r-\sigma}}=
\frac{\partial F}{\partial y}=F=0,x_{I_r-\sigma}\geq 0 \}. \]
We also call a function germ 
$F$ a generating family of a reticular Legendrian map
$\tilde{\pi} \circ i$ if $F$ is a generating family 
of the contact regular $r$-cubic configuration 
$\{ i(\tilde{L}^0_\sigma) \}_{ \sigma \subset I_r }$. 
}\end{dfn}

We call a contact diffeomorphism germ  $\phi$ 
on $(J^1({\mathbb R}^n,{\mathbb R}),0)$  
{\it a reticular diffeomorphism} if  
$\phi(\tilde{L}^0_\sigma)=\tilde{L}^0_\sigma$ for $\sigma \subset I_r$. 
We say that reticular Legendrian maps $\tilde{\pi} \circ i_1,\tilde{\pi}\circ i_2:(
\tilde{\mathbb L},0) \rightarrow (J^1({\mathbb R}^n,{\mathbb R}),0) 
\rightarrow ({\mathbb R}^n\times{\mathbb R},0)$ are 
{\it Legendrian equivalent} if there exist a reticular diffeomorphism $\phi$ and a Legendrian equivalence $\Theta$ of $\tilde{\pi}$ such that the following diagram is commutative:
\[
\begin{array}{ccccc}
(\tilde{\mathbb L},0) & \stackrel{i_1}{\longrightarrow} & (J^1({\mathbb R}^n,{\mathbb R}),0) & \stackrel{\tilde{\pi}}{\longrightarrow} & ({\mathbb R}^n\times{\mathbb R},0) \\
\phi|_{\tilde{\mathbb L}} \downarrow \ \ \ \ & & \Theta \downarrow \ \ & & g \downarrow \ \ \\
(\tilde{\mathbb L},0) & \stackrel{i_2}{\longrightarrow} & (J^1({\mathbb R}^n,{\mathbb R}),0) & \stackrel{\tilde{\pi}}{\longrightarrow} & ({\mathbb R}^n\times{\mathbb R},0), \\
\end{array}
\]
where $g$ is the diffeomorphism of the base space 
of $\tilde{\pi}$ induced from $\Theta$.

As after the definition of a reticular diffeomorphism in the previous part,
our definition of a reticular diffeomorphism and that in \cite{retLeg} 
are equivalent.\\

We say that function germs $F(x,y_1,\ldots,y_{k_1},q,z)
\in {\mathfrak M}(r;k_1+n+1)$ and $F(x,y_1,\ldots,y_{k_2},$ $q,z)\in 
{\mathfrak M}(r;k_2+n+1)$ are {\it stably  reticular 
${\cal P}$-${\cal K}$-equivalent} 
if $F$ and $G$ are reticular ${\cal P}$-${\cal K}$-equivalent 
after additions of non-degenerate quadratic forms in the variables $y$. 
\begin{tth}{\rm (cf., \cite[Theorem 5.6]{retLeg})}
{\rm (1)} For any reticular Legendrian map $\tilde{\pi} \circ i$, there exists a function germ $F\in {\mathfrak M}(r;k+n+1)$ which is a generating family of $\tilde{\pi}\circ i$. \\
{\rm (2)} For any C-non-degenerate function germ $F\in {\mathfrak M}(r;k+n+1)$, there exists a reticular Legendrian map of which $F$ is a generating family. \\
{\rm (3)} Two reticular Legendrian maps are Legendrian equivalent if and only if their generating families are stably reticular ${\cal P}$-${\cal K}$-equivalent. 
\end{tth}

\section{Stabilities of reticular Legendrian  maps}
\quad
We introduce several stabilities of reticular Legendrian maps:\\
{\bf Stability}:
For any open set $U$ in $J^1({\mathbb R}^n,{\mathbb R})$, 
we denote $C(U,J^1({\mathbb R}^n,{\mathbb R}))$ the space of contact embeddings from 
$U$ to $J^1({\mathbb R}^n,{\mathbb R})$ with $C^\infty$-topology. 
We say that 
a reticular Legendrian map $\tilde{\pi} \circ i:(\tilde{\mathbb L},0) \rightarrow 
(J^1({\mathbb R}^n,{\mathbb R}),0) \rightarrow({\mathbb R}^n\times{\mathbb R},0)$ is
 {\it stable} if the following
holds: For any extension $C$ of $i$ and any representative
$\tilde{C}\in C(U,J^1({\mathbb R}^n,{\mathbb R}))$ of $C$, there exists a neighbourhood
$N_{\tilde{C}}$ of $\tilde{C}$ such that 
for any $\tilde{C'}\in N_{\tilde{C}}$ 
the reticular Legendrian maps $\tilde{\pi} \circ \tilde{C}'_{w_0}|_{\tilde{\mathbb L}}$ and $\tilde{\pi} \circ i$ are Legendrian equivalent for some
$w_0=(0,\ldots,0,p^0_{r+1},\ldots,p^0_n)\in U$,
where the contact diffeomorphism germ $\tilde{C}'_{w_0}$ on 
$(J^1({\mathbb R}^n,{\mathbb R}),0)$ is chosen that 
the reticular Legendrian maps 
$\tilde{\pi} \circ \tilde{C}'|_{\tilde{\mathbb L}}:
(\tilde{\mathbb L},w_0) \rightarrow 
(J^1({\mathbb R}^n,{\mathbb R}),\tilde{C}'(w_0)) \rightarrow({\mathbb R}^n\times{\mathbb R},\tilde{\pi}\circ \tilde{C}'(w_0))$ 
and $\tilde{\pi} \circ \tilde{C}'_{w_0}|_{\tilde{\mathbb L}}:(\tilde{\mathbb L},0) \rightarrow 
(J^1({\mathbb R}^n,{\mathbb R}),0) \rightarrow(
{\mathbb R}^n\times{\mathbb R},0)$
are Legendrian equivalent.
\vspace{2mm}
\\
{\bf Homotopically stability}:
Let $\tilde{\pi} \circ i:(\tilde{\mathbb L},0) \rightarrow 
(J^1({\mathbb R}^n,{\mathbb R}),0)
\rightarrow({\mathbb R}^n\times{\mathbb R},0)$ be a reticular Legendrian map. 
A one-parameter family
of contact diffeomorphism germs  $\bar{C}:(J^1({\mathbb R}^n,{\mathbb R})
\times {\mathbb R},(0,0))\rightarrow
(J^1({\mathbb R}^n,{\mathbb R}),0)((Q,Z,P,t)\mapsto C_t(Q,Z,P))$
is called a {\it reticular Legendrian
  deformation} of $i$ if $C_0|_{\tilde{\mathbb L}}=i$.
Let $\phi$ be a reticular
diffeomorphism on $(J^1({\mathbb R}^n,{\mathbb R}),0)$. 
A map germ $\bar{\phi}:(J^1({\mathbb R}^n,{\mathbb R}) \times 
{\mathbb R},(0,0)) \rightarrow 
(J^1({\mathbb R}^n,{\mathbb R}),0)((Q,Z,P,t)\mapsto \phi_t(Q,Z,P))$
is called a {\it one-parameter
  deformation of reticular diffeomorphisms} of $\phi$ if $\phi_0=\phi$ and 
$\phi_t$ 
 is a reticular diffeomorphism for $t$ around $0$. 
We say that a reticular Legendrian map $\tilde{\pi} \circ i:({\mathbb
  L},0) \rightarrow (J^1({\mathbb R}^n,{\mathbb R}),0) \rightarrow
({\mathbb R}^n\times{\mathbb R},0)$ is {\it
  homotopically stable} if for any reticular Legendrian deformation 
$\bar{C}=\{ C_t\}$ of $i$,
there exists a one-parameter deformation of
reticular diffeomorphisms $\bar{\phi}=\{ \phi_t \}$ of 
$id_{J^1({\mathbb R}^n,{\mathbb R})}$
such that $C_t=\Theta_t\circ C_0\circ \phi_t$ for $t$ around
$0$.\vspace{2mm}
\\
{\bf Infinitesimally stability}:
Let $C$ be a contact diffeomorphism germ on $(J^1({\mathbb R}^n,{\mathbb R}),0)$.
We say that a vector field $v$ on $(J^1({\mathbb R}^n,{\mathbb R}),0)$ along $C$ 
is {\it infinitesimal contact transformation} if 
there exists a  reticular Legendrian deformation 
$\bar{C}=\{C_t\}$ on $(J^1({\mathbb R}^n,{\mathbb R}),0)$ 
such that $C_0=C$  and
$\frac{dC_t}{dt}|_{t=0}=v$.
We say that a vector field $\xi$ on $(J^1({\mathbb R}^n,{\mathbb R}),0)$ is 
{\it infinitesimally reticular diffeomorphism} if there exists 
a one-parameter deformation of
reticular diffeomorphisms  $\bar{\phi}=\{ \phi_t \}$ 
of $id_{J^1({\mathbb R}^n,{\mathbb R})}$ 
such that  $\frac{d\phi_t}{dt}|_{t=0}=\xi$.
We say that a vector field $\eta$ on  $(J^1({\mathbb R}^n,{\mathbb R}),0)$ is
{\it infinitesimally Legendrian equivalence} if 
there exists a
 one-parameter family of Legendrian equivalences $\bar{\Theta}=\{\Theta_t\}$ 
on $\tilde{\pi}$ such that $\Theta_0=id_{J^1({\mathbb R}^n,{\mathbb R})}$
and $\frac{d\Theta_t}{dt}|_{t=0}=\eta$.
We say that a reticular
Legendrian map $\tilde{\pi} \circ i:(\tilde{\mathbb L},0) \rightarrow 
(J^1({\mathbb R}^n,{\mathbb R}),0)
\rightarrow({\mathbb R}^n\times{\mathbb R},0)$ is  
 {\it infinitesimally stable} if for any
extension $C$ of $i$ and any infinitesimally contact transformation $v$ of $C$,
there exists an infinitesimally reticular diffeomorphism $\xi$ and 
an infinitesimally Legendrian equivalence $\eta$ such that $v=C_*\xi+\eta\circ C$.

\vspace{3mm}
We say that a function germ $H$ on $(J^1({\mathbb R}^n,{\mathbb R}),0)$ is {\it fiber preserving} if 
$H$ has the form $H(q,z,p)=\sum_{j=1}^nh_j(q,z)p_j+h_0(q,z)$.
\begin{tth}\label{stable:thleg}{\rm (cf., \cite[p.123 Theorem 7.4]{retLeg})}
Let $\tilde{\pi} \circ i:(\tilde{\mathbb L},0) \rightarrow 
(J^1({\mathbb R}^n,{\mathbb R}),0) \rightarrow({\mathbb R}^n\times{\mathbb R},0)$ 
be a reticular Legendrian map with a generating family $F(x,y,q,z)\in {\mathfrak M}(r;k+n+1)$.
 Then the following are equivalent.\\
{\rm (u)} $F$ is a reticular ${\cal P}$-${\cal K}$-stable unfolding of $F|_{q=z=0}$.\\
{\rm (hs)} $\tilde{\pi} \circ i$ is homotopically stable.\\
{\rm (is)} $\tilde{\pi} \circ i$ is infinitesimally stable.\\ 
{\rm (a)} For any function germ $f$ on $(J^1({\mathbb R}^n,{\mathbb R}),0)$, there exists a fiber preserving function germ $H$ on $ $ $(J^1({\mathbb R}^n,{\mathbb R}),0)$ such that $f\circ i=H\circ i$.\\
{\rm (s)} $\tilde{\pi} \circ i$ is stable.
\end{tth}
We shall prove the assertion (u)$\Leftrightarrow$(s) by a simpler way 
than that of \cite{retLeg} in Section \ref{GenLeg}.

\part{Generiticies of reticular Lagrangian, Legendrian maps}\label{Laggen:part}
\section{Finitely determinacy of reticular Lagrangian, Legendrian maps}
\begin{dfn}{\rm
Let $\pi\circ i:({\mathbb L},0) \rightarrow (T^* {\mathbb R}^n,0) 
\rightarrow({\mathbb R}^n,0)$ be a reticular Lagrangian map
and $l$ be a non-negative number.
We say that $\pi\circ i$ is $l$-determined if the following condition holds:
For any extension $S$ of $i$, the reticular Lagrangian map 
$\pi\circ S'|_{{\mathbb L}}$ and $\pi\circ i$ are Lagrangian 
equivalent for any symplectic diffeomorphism germ $S'$ on 
$(T^* {\mathbb R}^n,0)$ satisfying $j^lS(0)=j^lS'(0)$.
}\end{dfn}
\begin{lem}\label{SS:lem}
Let $\pi\circ i:({\mathbb L},0) \rightarrow (T^* {\mathbb R}^n,0) 
\rightarrow({\mathbb R}^n,0)$ be a reticular Lagrangian map.
Let $S_1,S_2$ are symplectic diffeomorphism germ on 
$(T^* {\mathbb R}^n,0)$ such that 
$S_1|_{{\mathbb L}}=S_2|_{{\mathbb L}}=i$.
Then there exists a reticular diffeomorphism germ $\phi$ such that 
$S_1=S_2\circ \phi$.
\end{lem}
{\it Proof.}
Set $\phi=(S_2)^{-1}\circ S_1$.
Then we have that $\phi|_{{\mathbb L}}=id$ and $S_1=S_2\circ \phi$\hfill $\blacksquare$

By this lemma we have that the finitely determinacy of reticular Lagrangian maps 
do not depend on choices of extensions of their reticular Lagrangian embeddings.
\begin{tth}\label{n+1det:Lag}
Let $\pi\circ i:({\mathbb L},0) \rightarrow (T^* {\mathbb R}^n,0) 
\rightarrow({\mathbb R}^n,0)$ be a reticular Lagrangian map.
If $\pi\circ i$ is infinitesimally stable then $\pi\circ i$ is $(n+1)$-determined.
\end{tth}
{\it Proof.}
Let $S$ be an extension of $i$.
Since the infinitesimally stability of reticular Lagrangian maps is invariant under  
Lagrangian equivalences,
we may assume that the canonical relation $P_S$ associated with $S$ 
has the form:
\[ P_S=\{(Q,-\frac{\partial H}{\partial Q}(Q,p),
-\frac{\partial H}{\partial p}(Q,p),p)\in 
(T^* {\mathbb R}^n\times T^* {\mathbb R}^n,(0,0))
\}\]
for some function germ $H(Q,p)\in {\mathfrak M}(2n)^2$.
Then $F(x,y,q)=H_0(x,y)+\langle y,q\rangle\in {\mathfrak M}(r;n+n)^2$
is a generating family of $\pi\circ i$,
where $H_0(x,y)=H(x_1,\ldots,x_r,0,\ldots,0,y_1,\ldots,y_n)\in {\mathfrak M}(r;n)^2$.
Since $\pi\circ i$ is infinitesimally stable, it follows that $F$ is a
reticular ${\cal P}$-${\cal R}^+$-infinitesimally versal unfolding of 
$H_0(x,y)$ by Theorem \ref{stable:thlag}.
This means that ${\cal E}(r;n)/\langle x\frac{\partial H_0}{\partial x},
\frac{\partial H_0}{\partial y}\rangle_{{\cal E}(r;n)}$ is at most $(n+1)$-dimension.
It follows that
\[{\mathfrak M}(r;n)^{n+1}\subset \langle x\frac{\partial H_0}{\partial x},
\frac{\partial H_0}{\partial y}\rangle_{{\cal E}(r;n)}.\]
Therefore we have that
\[ {\mathfrak M}(r;n)^{n+3}\subset {\mathfrak M}(r;n)(
\langle x\frac{\partial H_0}{\partial x}\rangle +{\mathfrak M}(r;n)
\langle \frac{\partial H_0}{\partial y}\rangle).\]
This means that $H_0$ is reticular ${\cal R}$-$(n+2)$-determined
by Lemma \ref{findet:lm}.
Let a symplectic diffeomorphism germ $S'$ on $(T^* {\mathbb R}^n,0)$ 
satisfying $j^{n+1}S(0)=j^{n+1}S'(0)$ be given.
Since $\frac{\partial p\circ S}{\partial p}(0)=
\frac{\partial p\circ S'}{\partial p}(0)$, it follows that 
there exists a function germ $H'(Q,p)\in {\mathfrak M}(2n)^2$ such that 
\[ P_{S'}=\{(Q,-\frac{\partial H'}{\partial Q}(Q,p),
-\frac{\partial H'}{\partial p}(Q,p),p)\}.\]
Then the function germ $G(x,y,q):=H'_0(x,y)+\langle y,q\rangle\in 
{\mathfrak M}(r;n+n)^2$
is a generating family of $\pi\circ S'|_{{\mathbb L}}$,
where $H_0'(x,y)=H'(x,0,y)\in {\mathfrak M}(r;n)^2$.
Since $j^{n+1}S(0)=j^{n+1}S'(0)$, it means that 
$j^{n+1}\frac{\partial H}{\partial Q}(0)=
j^{n+1}\frac{\partial H'}{\partial Q}(0),\ 
j^{n+1}\frac{\partial H}{\partial p}(0)=
j^{n+1}\frac{\partial H'}{\partial p}(0)$.
Therefore we have that $j^{n+2}H(0)=j^{n+2}H'(0)$
and hence $j^{n+2}H_0(0)=j^{n+2}H_0'(0)$.
Thus there exists $\Phi(x,y)\in {\cal B}(r;n)$ such that 
$H_0=H_0'\circ \Phi$.
We set $G'(x,y,q):=G(\Phi(x,y),q)\in {\mathfrak M}(r;n+n)^2$.
Then $G$ and $G'$ are reticular ${\cal P}$-${\cal R}$-equivalent,
and 
$F$ and $G'$ are reticular ${\cal P}$-${\cal R}^+$-infinitesimal versal 
unfoldings of $H_0(x,y)$.
It follows that $F$ and $G$ are reticular ${\cal P}$-${\cal R}^+$-equivalent.
Therefore $\pi\circ i$ and $\pi \circ S'|_{{\mathbb L}}$ are
Lagrangian equivalent.\hfill $\blacksquare$

\vspace{5mm}
We  consider contact diffeomorphism 
germs on $(J^1({\mathbb R}^n,{\mathbb R}),0)$.
 Let $(Q,Z,P)$ be canonical coordinates on the 
source space
and $(q,z,p)$ be canonical coordinates of the target space. We define the 
following notations:\\
$\imath:(J^1({\mathbb R}^n,{\mathbb R})\cap \{Z=0 \},0)\rightarrow 
(J^1({\mathbb R}^n,{\mathbb R}),0)$ be the inclusion map on the source space.
\begin{eqnarray*}
  C(J^1({\mathbb R}^n,{\mathbb R}),0) & = & \{ C:(J^1({\mathbb R}^n,{\mathbb R}),0)
\rightarrow (J^1({\mathbb R}^n,{\mathbb R}),0) |C:\ {\rm contact\ diffeomorphism }
 \} \\
 C^\theta (J^1({\mathbb R}^n,{\mathbb R}),0) & = & \{ C\in C(J^1({\mathbb R}^n
,{\mathbb R}),0)|
\mbox{ $C$ preserves the canonical $1$-form }\ \} \\
 C_Z (J^1({\mathbb R}^n,{\mathbb R}),0) & = & \{ C\circ\imath\ |C \in 
C(J^1({\mathbb R}^n,{\mathbb R}),0) \}\\
 C^\theta_Z (J^1({\mathbb R}^n,{\mathbb R}),0) & = & \{ C\circ\imath\ |
C \in C^\theta (J^1({\mathbb R}^n,{\mathbb R}),0) \}.
 \end{eqnarray*}
\begin{dfn}\label{l,l+1detLeg}{\rm
Let $\tilde{\pi}\circ i:(\tilde{\mathbb L},0) \rightarrow (J^1({\mathbb R}^n,{\mathbb R}),0) 
\rightarrow({\mathbb R}^n\times {\mathbb R},0)$ be a reticular Legendrian map.
We say that $\tilde{\pi}\circ i$ is $l$-determined if the following condition holds:
For any extension $C\in C(J^1({\mathbb R}^n,{\mathbb R}),0)$ of $i$, the reticular Legendrian map 
$\tilde{\pi}\circ C'|_{\tilde{\mathbb L}}$ and $\tilde{\pi}\circ i$ are  Legendrian 
equivalent for $C'\in C(J^1({\mathbb R}^n,{\mathbb R}),0)$ satisfying $j^lC(0)=j^lC'(0)$.
}\end{dfn}

\begin{lem}{\rm (cf., \cite[p.116 Lemma 7.2]{retLeg})}\label{exleg:lm}
Let $\{L_\sigma \}_{\sigma \subset I_r}$ be a contact regular $r$-cubic configuration in 
$(J^1({\mathbb R}^n,{\mathbb R}),0)$ defined by $C\in C(J^1({\mathbb R}^n,
{\mathbb R}),0)$. Then there exists $C'\in C^\theta (J^1({\mathbb R}^n,$ ${\mathbb R}),0)$
 that also defines $\{L_\sigma \}_{\sigma \subset I_r}$.
\end{lem}

By Lemma \ref{exleg:lm}  we may consider the following other definition
of finitely determinacy of reticular Legendrian maps:\\
(1) The definition given by replacing $C(J^1({\mathbb R}^n,{\mathbb R}),0)$
to $C^\theta(J^1({\mathbb R}^n,{\mathbb R}),0)$.\\
(2) The definition given by replacing $C(J^1({\mathbb R}^n,{\mathbb R}),0)$
to $C_Z(J^1({\mathbb R}^n,{\mathbb R}),0)$.\\
(3) The definition given by replacing $C(J^1({\mathbb R}^n,{\mathbb R}),0)$
to $C^\theta_Z(J^1({\mathbb R}^n,{\mathbb R}),0)$.\\
Then the following holds:
\begin{prop}
Let $\tilde{\pi}\circ i:(\tilde{\mathbb L},0) \rightarrow (J^1({\mathbb R}^n,{\mathbb R}),0)
\rightarrow({\mathbb R}^n\times {\mathbb R},0)$ be a reticular Legendrian map.
Then \\
{\rm (A)} If $\tilde{\pi}\circ i$ is $l$-determined of the original definition, then $\tilde{\pi}\circ i$ is $l$-determined of the 
definition {\rm (1)}.\\
{\rm (B)} If $\tilde{\pi}\circ i$ is $l$-determined of the definition {\rm (1)}, then $\tilde{\pi}\circ i$ is $l$-determined of the definition {\rm (3)}.\\
{\rm (C)} If $\tilde{\pi}\circ i$ is $l$-determined of the definition {\rm (3)}, then $\tilde{\pi}\circ i$ is $(l+1)$-determined of the definition {\rm (2)}.\\
{\rm (D)} If $\tilde{\pi}\circ i$ is $l$-determined of the definition {\rm (2)}, then $\tilde{\pi}\circ i$ is $l$-determined of the original definition.
\end{prop}
{\it Proof.}
(A) Let $C\in C^\theta(J^1({\mathbb R}^n,{\mathbb R}),0)$ be an 
extension of $i$. 
Let $C'\in C^\theta(J^1({\mathbb R}^n,{\mathbb R}),0)$ satisfying 
$j^lC(0)=j^lC'(0)$ be given.
Since $C,C'\in C(J^1({\mathbb R}^n,{\mathbb R}),0)$,
we have that $\tilde{\pi}\circ C'|_{\tilde{{\mathbb L}}}$ and $\tilde{\pi}\circ i$ are  Legendrian 
equivalent.\\
(B) Let $C\in C^\theta_Z(J^1({\mathbb R}^n,{\mathbb R}),0)$ be an 
extension of $i$. 
Let $C'\in C^\theta_Z(J^1({\mathbb R}^n,{\mathbb R}),0)$ satisfying 
$j^lC(0)=j^lC'(0)$ be given.
We define $C_1,C'_1\in C^\theta(J^1({\mathbb R}^n,{\mathbb R}),0)$ by 
 $C_1(Q,Z,P):=(q_C(Q,P),Z+z_C(Q,P),p_C(Q,P)),
C'_1(Q,Z,P):=(q_{C'}(Q,P),Z+z_{C'}(Q,P),p_{C'}(Q,P))$.
Then $C_1$ is an extension of $i$ and we have that $j^lC_1(0)=j^lC_1'(0)$.
Therefore we have that 
$\tilde{\pi}\circ C_1'|_{\tilde{{\mathbb L}}}=
\tilde{\pi}\circ C'|_{{\tilde{\mathbb L}}}$ and $\tilde{\pi}\circ i$ are  Legendrian 
equivalent.\\
(D) Let $C\in C(J^1({\mathbb R}^n,{\mathbb R}),0)$ be an 
extension of $i$. 
Let $C'\in C(J^1({\mathbb R}^n,{\mathbb R}),0)$ satisfying 
$j^lC(0)=j^lC'(0)$ be given.
We set $C_1:=C|_{Z=0},C'_1:=C'|_{Z=0}\in C_Z(J^1({\mathbb R}^n,{\mathbb R}),0)$
and we have that $j^lC_1(0)=j^lC_1'(0)$.
Therefore $\tilde{\pi}\circ C_1'|_{\tilde{{\mathbb L}}}=
\tilde{\pi}\circ C'|_{\tilde{{\mathbb L}}}$ and $\tilde{\pi}\circ i$ are  Legendrian 
equivalent.\\
(C) Let $C\in C_Z(J^1({\mathbb R}^n,{\mathbb R}),0)$ be an 
extension of $i$. 
Let $C'\in C_Z(J^1({\mathbb R}^n,{\mathbb R}),0)$ satisfying 
$j^{l+1}C(0)=j^{l+1}C'(0)$ be given.
Then there exist function germs $f(Q,P),g(Q,P)\in {\cal E}(2n)$
such that $C^{*}(dz-pdq)=-fPdQ,C^{'*}(dz-pdq)=-gPdQ$.
Indeed $f$ is defined by that 
$fP_j=-\frac{\partial z_C}{\partial Q_j}+p_C
\frac{\partial q_C}{\partial Q_j}$ for $j=1,\ldots,n$.
We define the diffeomorphism germs $\phi,\psi$ on 
$(J^1({\mathbb R}^n,{\mathbb R})\cap \{Z=0\},0)$ by
$\phi(Q,P)=(Q,fP),\psi(Q,P)=(Q,gP)$.
We set $C_1:=C\circ \phi^{-1},C'_1:=C'\circ\psi^{-1}\in 
C^\theta_Z(J^1({\mathbb R}^n,{\mathbb R}),0)$
Then $j^l\phi(0)$ and $j^l\psi(0)$ depend only on $j^{l+1}C(0)$,
therefore we have that $j^lC_1(0)=j^lC_1'(0)$.
Since $\tilde{\pi}\circ i$ and $\tilde{\pi}\circ C_1|_{\tilde{{\mathbb L}}}$ are Legendrian 
equivalent,
it follows that 
$\tilde{\pi}\circ C_1|_{\tilde{{\mathbb L}}}$ and 
$\tilde{\pi}\circ C'_1|_{\tilde{{\mathbb L}}}$ are  Legendrian 
equivalent.
Therefore we have that $\tilde{\pi}\circ i$ and
 $\tilde{\pi}\circ C'|_{\tilde{{\mathbb L}}}$ are Legendrian 
equivalent.\hfill $\blacksquare$
\begin{tth}\label{n+1det:Leg}
Let $\tilde{\pi}\circ i:(\tilde{{\mathbb L}},0) \rightarrow 
(J^1({\mathbb R}^n,{\mathbb R}),0) 
\rightarrow({\mathbb R}^n\times {\mathbb R},0)$ be a reticular Legendrian map.
If $\tilde{\pi}\circ i$ is infinitesimally stable then $\tilde{\pi}\circ i$ is $(n+3)$-determined.
\end{tth}
{\it Proof.}
It is enough to prove $\tilde{\pi}\circ i$ is $(n+2)$-determined of Definition 
\ref{l,l+1detLeg} (3) and 
this is proved by an analogous method of Theorem \ref{n+1det:Lag}.
We write a sketch of the proof.
Let $C\in C^\theta_Z(J^1({\mathbb R}^n,{\mathbb R}),0)$ be an 
extension of $i$.
Then we may assume that $P_C$ has the form 
\[ P_C=\{(Q,-\frac{\partial H}{\partial Q}(Q,p);
H(Q,p)-\langle \frac{\partial H}{\partial Q}(Q,p),p\rangle,
-\frac{\partial H}{\partial p}(Q,p),p)\}\]
for some function germ $H(Q,p)\in {\mathfrak M}(2n)^2$.
Then $F(x,y,q,z)=H_0(x,y)+\langle y,q\rangle -z \in {\mathfrak M}(r;n+n+1)$
is a generating family of $\tilde{\pi}\circ i$,
where $H_0(x,y)=H(x,0,y)\in {\mathfrak M}(r;n)^2$.
Then $F$ is a
reticular ${\cal P}$-${\cal K}$-infinitesimally stable unfolding of 
$H_0(x,y)$.
It follows that $H_0$ is reticular ${\cal K}$-$(n+2)$-determined.
Let $C'\in C^\theta_Z(J^1({\mathbb R}^n,{\mathbb R}),0)$ 
satisfying $j^{n+2}C(0)=j^{n+2}C'(0)$ be given.
There exists a function germ $H'(Q,p)\in {\mathfrak M}(2n)$ such that 
\[ P_{C'}=\{(Q,-\frac{\partial H'}{\partial Q}(Q,p);
H'(Q,p)-\langle \frac{\partial H'}{\partial Q}(Q,p),p\rangle,
-\frac{\partial H'}{\partial p}(Q,p),p)\}.\]
Since $H=z-qp$ on $P_C$ and $H'=z-qp$ on $P_{C'}$,
we have that $j^{n+2}H_0(0)=j^{n+2}H_0'(0)$, where 
$H_0'(x,y)=H'(x,0,y)\in {\mathfrak M}(r;n)^2$.
The function germ $G(x,y,q,z)=H_0'(x,y)+\langle y,q\rangle -z
\in {\mathfrak M}(r;n+n+1)$
is a generating family of $\tilde{\pi}\circ C'|_{\tilde{{\mathbb L}}}$.
Then there exist $\Phi(x,y)\in {\cal B}(r;n)$ and an unit $a\in
{\cal E}(r;n)$ such that 
$H_0=a\cdot H_0'\circ \Phi$.
We set $G'(x,y,q,z)=a(x,y)G(\Phi(x,y),q,z)\in {\mathfrak M}(r;n+n+1)$.
Then $G$ and $G'$ are reticular ${\cal P}$-${\cal K}$-equivalent
and $F$ and $G'$ are reticular ${\cal P}$-${\cal K}$-infinitesimal stable 
unfoldings of $H_0(x,y)$.
It follows that $F$ and $G$ are reticular ${\cal P}$-${\cal K}$-equivalent.
Therefore $\tilde{\pi}\circ i$ and $\tilde{\pi} \circ C'|_{\tilde{{\mathbb L}}}$ are
Legendrian equivalent.\hfill $\blacksquare$

\vspace{3mm}
In order to prove Theorem \ref{transv} and Theorem \ref{transv:leg},
we require the following lemmas:
\begin{lem}{\rm (cf., \cite[p.53 Lemma 4.6]{golubitsky:text})}
\label{translem:lem}
Let $V,W$  be smooth manifolds with $Q$ a submanifold of $W$,
$F$ be a topological space.
We equip $C^\infty (V,W)$ with the Whitney $C^\infty$-topology. 
Let $j:F\rightarrow  C^\infty (V,W)$ be a map
{\rm  (}not necessary continuous\rm{)}.
We suppose that: For each $f\in F$,
there exist a manifold $E$, $e_0\in E$, and a continuous map
$\phi:E\rightarrow F,\ \phi(e_0)=f$ such that
the map $\Phi:E\times V\rightarrow W,\ 
\Phi(e,x)=j( \phi(e))(x)$, is smooth and 
transversal to $Q$.
Then the set 
\[ T=\{ f\in F\ | \ j(f)\mbox{ is transversal to } Q\} \]
is dense in $F$.
\end{lem}
\begin{lem}\label{transopen:lem}
Let $V,W$  be smooth manifolds with $Q$ a submanifold of $W$,
and $K\subset Q$ be a compact set.
Then the set 
\[ T_K=\{ f\in C^\infty(V,W)\ | \ j^lf\mbox{ is transversal to } Q
\mbox{ on } K\} \]
is open in $C^\infty(V,W)$ with Whitney $C^{l+1}$-topology 
\rm{(}{\it hence Whitney $C^\infty$-topology}\rm{)}.
\end{lem}

\section{Genericity of reticular Lagrangian maps}\label{GeneLag}
\quad
Let  $J^l(2n,2n)$ be the set of $l$-jets of map germs  from 
$(T^*{\mathbb R}^n,0)$ 
to $(T^*{\mathbb R}^n,0)$ and 
$S^l(n)$ be the Lie group in $J^l(2n,2n)$ 
consists of $l$-jets of symplectic diffeomorphism germs 
on $(T^*{\mathbb R}^n,0)$.

We consider the Lie group $L^l(2n)\times L^l(2n)$ acts on $J^l(2n,2n)$ as coordinate changes of 
the source and target spaces.
We also consider the Lie subgroup $rLa^l(n)$ of $L^l(2n)\times L^l(2n)$ consists of 
 $l$-jets of reticular 
diffeomorphisms on the source space and 
$l$-jets of Lagrangian equivalences of $\pi$: 
\begin{eqnarray*}
rLa^l(n)=\{ (j^l\phi(0),j^l\Theta(0))\in L^l(2n)\times L^l(2n)\ | \ \phi
 \mbox{ is a reticular diffeomorphism on } \\
(T^*{\mathbb R}^n,0), \Theta \mbox{ is a Lagrangian equivalence of } 
\pi \}. 
\end{eqnarray*}
The group $rLa^l(n)$ acts on $J^l(2n,2n)$ and  $S^l(n)$
is invariant under this action.
Let $S$ be  a symplectic diffeomorphism germ on $(T^*{\mathbb R}^n,0)$ and set 
$z=j^lS(0)$.
We denote  the orbit $rLa^l(n)\cdot z$ by $[z]$.
Then 
\[ [z]=\{ j^lS'(0)\in S^l(n)\ | \ \pi \circ i \mbox{ and } \pi\circ S'|_{\mathbb L} 
\mbox{ are Lagrangian equivalent} \}. \]

In this section we denote by $X_f$ the Hamiltonian vector field on 
$(T^*{\mathbb R}^n,0)$ for a function germ $f$ on $(T^*{\mathbb R}^n,0)$.
That is
\[ X_f=\sum_{j=1}^n(\frac{\partial f}{\partial p_j}\frac{\partial }{\partial q_j}
-\frac{\partial f}{\partial q_j}\frac{\partial }{\partial p_j}).\]

We denote by $VI_S$  the vector space consists of infinitesimal symplectic transformations of $S$
and  denote by $VI^0_S$ the subspace of $VI_S$ consists of germs which vanishes on $0$.
We denote by $VL_{T^*{\mathbb R}^n}$ by the vector space consists of infinitesimal Lagrangian 
equivalences of $\pi$  and denote by $VL^0_{T^*{\mathbb R}^n}$
by the subspace of $VL_{T^*{\mathbb R}^n}$ consists of germs which vanishes at $0$.

We denote by $V^0_{\mathbb L}$ the vector space consists of infinitesimal reticular diffeomorphisms 
on $(T^*{\mathbb R}^n,0)$ which vanishes at $0$:
\[ V^0_{\mathbb L}=\{ \xi \in X(T^*{\mathbb R}^n,0)\ | \ 
\xi \mbox{ is tangent to } L^0_\sigma \mbox { for all } \sigma\subset I_r,\ \xi(0)=0\}. \]

\vspace{5mm}
From now on, we denote by ${\cal E}_{T^* {\mathbb R}^n}$  the ring of 
smooth function germs on $(T^* {\mathbb R}^n,0)$ and denote 
by ${\mathfrak M}_{T^* {\mathbb R}^n}$ its maximal ideal. 
We also denote other notations analogously.
\begin{lem}\label{infsta:Laglem}
{\rm (1)} A vector field germ $v$ on $(T^*{\mathbb R}^n,0)$ along $S$ is 
belongs to $VI_S$ if and only if there exists a function germ $f$ on $(T^*{\mathbb R}^n,0)$
such that $v=X_f\circ S$,\\
{\rm (2)} A vector field germ $\eta$ on $(T^*{\mathbb R}^n,0)$ is
belongs to $VL_{T^*{\mathbb R}^n}$ if and only if there exists a fiber preserving function germ
$H$ on $(T^*{\mathbb R}^n,0)$ such that $\eta=X_H$.\\
{\rm (3)}  A vector field germ $\xi$ on $(T^*{\mathbb R}^n,0)$ is
belongs to $V^0_{\mathbb L}$ if and only if there exists a function germ $g\in B_0$
such that $\xi=X_g$, where $B_0=\langle q_1p_1,\ldots,q_rp_r\rangle_{{\cal E}_{T^* {\mathbb R}^n}}+
{\mathfrak M}_{T^* {\mathbb R}^n}\langle q_{r+1},\ldots,q_n\rangle$
is a submodule of ${\cal E}_{T^* {\mathbb R}^n}$.
\end{lem}
By this lemma we have that:
\[ VI^0_S=\{ v:(T^*{\mathbb R}^n,0)\rightarrow (T(T^*{\mathbb R}^n),0)\ | \
v=X_f\circ S \mbox{ for some }f\in {\mathfrak M}^2_{T^*{\mathbb R}^n}\},\]
\[ VL^0_{T^*{\mathbb R}^n}=\{ \eta\in X(T^*{\mathbb R}^n,0)\ | \
\eta=X_H \mbox{ for some fiver preserving function germ } H\in {\mathfrak M}^2_{T^*{\mathbb R}^n}\}, \]
\[ V^0_{\mathbb L}=\{ \xi \in X(T^*{\mathbb R}^n,0)\ | \ 
\xi=X_g  \mbox{ for some }g\in B_0\}. \]

\vspace{2mm}
We define the homomorphism $tS:VI^0_{\mathbb L}\rightarrow 
VI^0_S$ by $tS(v)=S_*v$ and define the homomorphism 
$wS:VL^0_{T^*{\mathbb R}^n}\rightarrow VI^0_S$ by 
$wS(\eta)=\eta\circ S$.
\begin{lem}
Let $S$ be  a symplectic diffeomorphism germ on $(T^*{\mathbb R}^n,0)$ and set 
$z=j^lS(0)$. Then
\[ T_z(rLa^l(n)\cdot z)=\pi_l( tS( V^0_{\mathbb L})+wS(VL^0_{T^*{\mathbb R}^n})).\]
\end{lem}

We denote $VI^{l}_S$ the subspace of $VI_S$ consists of infinitesimal
symplectic transformation germs of $S$ whose $l$-jets are $0$:
\[  VI^l_S=\{ v\in VI_S\ | \ j^lv(0)=0\}. \]
We consider the surjective projection $\pi_l:VI_S\rightarrow T_z(S^l(n))$. 
Since $(j^lS)_*(\frac{ \partial }{\partial q_j})=\pi_l( S_*\frac{\partial }{\partial q_j}),$ $
(j^lS)_*(\frac{ \partial }{\partial p_j})=\pi_l( S_*\frac{\partial }{\partial p_j})$,
it follows that 
$j^lS$ is transversal to $[z]$ if and only if 
$(j^lS)_*(T_0(T^*{\mathbb R}^n))$ $+T_z[z]=T_z(S^l(n))$ and this holds 
if and only if 
\[ (\pi_l)^{-1}((j^lS)_*(T_0(T^*{\mathbb R}^n))+
tS( V^0_{\mathbb L})+wS(VL^0_{T^*{\mathbb R}^n}))+VI^{l+1}_S=VI_S
\] 
 and 
this holds if and only if 
\[ 
tS( V^0_{\mathbb L})+wS(VL_{T^*{\mathbb R}^n}))+VI^{l+1}_S=VI_S. \]

\vspace{3mm}
Let $N,M$ be $2n$ dimensional symplectic manifolds.
We denote $S(N,M)$ the space of symplectic embeddings from
$N$ to $M$ with the topology induced from the Whitney $C^\infty$-topology
of $C^\infty(N,M)$.
We define that
\begin{eqnarray*}
 J^l_S(N,M)=\{ j^lS(u_0)\in J^l(N,M)\ | \hspace{9cm} \\
 S:(N,u_0)\rightarrow M \mbox{ 
is a symplectic embedding germ},\ u_0\in N\}.
\end{eqnarray*}
 \begin{prop}
$S(N,M)$ is a Baire space.
\end{prop}
This is proved  by an analogous method
of the assertion that $C^\infty(N,M)$ is a Baire space
(cf., \cite[p.44 Proposition 3.3]{golubitsky:text}).
\begin{tth}[Symplectic transversality theorem]\label{symplectictrans:tth}
Let $N,M$ be $2n$ dimensional symplectic manifolds.
 Let 
$Q_j, j=1,2,\ldots$ be submanifolds of $J^l_S(N,M)$.
Then the set
\[ T=\{ S\in S(N,M)\ | \ j^lS \mbox{ is transversal to } Q_j 
\mbox{ for all } j\in {\mathbb N}\} \]
is residual set in $S(N,M)$. In particular $T$ is dense.
\end{tth}
{\it Proof}.
We apply for Lemma \ref{translem:lem} that 
$V=N,W=J^l_S(N,M)$, and $F=S(N,M)$.
We reduce our assertion to local situations by 
choosing a countable covering of $Q_j$ by 
 sufficiently small compact sets $K_{j,k}$'s. 
Then the sets $T_{j,k}=\{ S\in S(N,M)\ | \ j^lS \mbox{ is transversal to } Q_j 
\mbox{ at } K_{j,k} \}$ are open set by Lemma \ref{translem:lem} and 
we have that $T=\cap T_{j,k}$
We fix a symplectic embedding $S\in S(N,M)$.
For each $u_0\in N$ there exist local symplectic coordinate 
systems  of 
$N$ around $u_0$ and $M$ around $S(u_0)$ such that 
$S$  is given by $(q,p)\mapsto (q,p)$ around $0$.

For each $j,k$ we take $E$ by a sufficiently small neighborhood of $0$
in $P(2n,1;l+1)$ and take a smooth function $\rho:T^*{\mathbb R}^n
\rightarrow [0,1]$ such that $\rho$ 
is $1$ on a neighborhood of $0$ and 
zero outside a 
compact set,
where $P(2n,1;l+1)$ is the set of not higher than $(l+1)$-degree 
polynomials on $2n$ variables.

For each $H\in E$ we define
$H'(Q,p)=\rho(Q,p)H(Q,p)-\langle Q,p\rangle$
and  
 $\psi_{H'}(Q,p)=
(Q,$ $-\frac{\partial H'}{\partial Q}(Q,p))$
for $(Q,p)\in T^*{\mathbb R}^n$ around $0$.
Then 
there exists a neighbourhood $U$ of $0$ in $T^*{\mathbb R}^n$
such that 
$\psi_{H'}$ is a embedding on $U$ 
and equal to the identity map  outside a compact set for any $H\in E$.
Therefore 
there exists a neighbourhood $U'$ of $0$ in $T^*{\mathbb R}^n$ 
such that the map 
$E\rightarrow C^\infty(U',U),\ H\mapsto (\psi_H^{-1})|_{U'}$ 
is well defined and continuous. Each $(\psi_H^{-1})|_{U'}$ is 
equal the identity map outside a compact set around $0$.
We set  that 
\[ \phi(H)(Q,P)=(
-\frac{\partial H'}{\partial p}(Q,p),p)\circ 
(\psi_{H'})^{-1}(Q,P)\mbox{ for } (Q,P)\in U'.\]
Then $\phi(H)$ is a symplectic diffeomorphism around $u_0$ 
which has the canonical relation with the generating function $H'(Q,p)$
and equal to $S$  outside a compact set.
It follows that 
the source space of 
$\phi(H)$ may be extended naturally to $E\times N$.
We also denote this by $\phi_H\in S(N,M)$.
Then the map
\[ \Phi:E\times N\rightarrow J^l_S(N,M), \Phi(H,q,p)=j^l(\phi_H)(q,p)\]
is a submersion around $(0,u_0)$
Therefore $\Phi$ is transversal to $K_{j,k}$.
So we have the result.\hfill $\blacksquare$

\vspace{3mm}
For $\tilde{S}\in S(U,T^*{\mathbb R}^n)$ we define
 the continuous map 
$j^l_0\tilde{S}:U\rightarrow S^l(n)$ by $w$ to the $l$-jet of 
$\tilde{S}_w$ at $0$.
We remark that $j^l\tilde{S}_w(0)$ has the form 
$j^l\tilde{S}_w(0)=(w,\tilde{S}(w),j^l_0\tilde{S}_w(0))$.

\begin{tth}\label{transv}
Let $U$ be a neighborhood of $0$ in $T^*{\mathbb R}^n$, and let
 $Q_1,Q_2,\ldots$ are submanifolds of $S^l(n)$. Then the set
$L_U=\{ (q,p)\in U |q=p_1=\cdots =p_r=0 \}$ and 
\[ T=\{ \tilde{S}\in S(U,T^*{\mathbb R}^n)| j^l_0\tilde{S}
 \mbox{ is transversal to } Q_j \mbox{ on } L_U
\mbox{ for all } j \} \]
is a residual set in $S(U,T^*{\mathbb R}^n)$.
\end{tth}
{\it Proof}. 
We set  $Q'_j=L_U\times T^* {\mathbb R}^n\times Q_j\subset
J^l_S(U,T^* {\mathbb R}^n)$.
We choose a countable covering of $Q'_j$ by  sufficiently small compact 
sets $K_{j,k}{}'s$  for all $j$.
We apply $N=V,M=T^*{\mathbb R}^n$ for Theorem \ref{symplectictrans:tth}.
We have that 
\[  T=\{ \tilde{S}
\in S(U,T^*{\mathbb R}^n)| j^l\tilde{S} \mbox{ is transversal to } 
Q'_j \mbox{ on } K_{j,k} \mbox{ for all } j,k \}. \]
It follows that $T$ is a residual set in $S(U,T^*{\mathbb R}^n)$.
\hfill $\blacksquare$
\begin{tth}\label{ts:thlag}
Let $\pi\circ i:({\mathbb L},0) \rightarrow (T^* {\mathbb R}^n,0) 
\rightarrow({\mathbb R}^n,0)$ be a reticular Lagrangian map. 
Let $S$ be an extension of $i$ and $l\geq n+1$.
Let $B=\langle q_1p_1,\ldots,q_rp_r, q_{r+1},\ldots,q_n\rangle_{{\cal E}_{T^* {\mathbb R}^n}}$
be a submodule of ${\cal E}_{T^* {\mathbb R}^n}$.
Then the followings are equivalent:\\
{\rm (s)} $\pi\circ i$ is stable.\\
{\rm (t)} $j^l_0S$ is transversal to $[j^l_0S(0)]$ at $0$.\\
{\rm (a')} ${\cal E}_{T^* {\mathbb R}^n}/(B+{\mathfrak M}_{T^* {\mathbb R}^n}^{l+2})$ is generated by 
$1,p_1\circ S,\ldots,p_n\circ S$ as ${\cal E}_ {q}$-module via $\pi\circ S$.\\
{\rm (a'')} ${\cal E}_{T^* {\mathbb R}^n}/B$ is generated by 
$1,p_1\circ S,\ldots,p_n\circ S$ as ${\cal E}_ {q}$-module via $\pi\circ S$.\\
{\rm (is)} $\pi\circ i$ is infinitesimally stable. 
\end{tth}
{\ Proof}.
(s)$\Rightarrow$(t). Let $\tilde{S}\in S(U,T^* {\mathbb R}^n)$ be a representative 
of $S$.
By theorem \ref{transv} there exists a symplectic embedding $\tilde{S'}$ 
around $\tilde{S}$ such that 
$j^l_0{\tilde{S'}}$ is transversal to $[j^l_0S(0)]$ at 
$w=(0,\ldots,0,p_{r+1},\ldots,p_n)\in U$.
Since $\pi \circ i$ is stable, $\pi \circ i$ and $\pi \circ \tilde{S'_w}|_{\mathbb L}$ are
Lagrangian equivalent.
This means that $[j^l_0\tilde{S'_w}(0)]=[j^l_0S(0)]$ and hence 
$j^l_0S$ is transversal to $[j^l_0S(0)]$ at $0$.\\
(t)$\Leftrightarrow$(a').
By Lemma \ref{infsta:Laglem}, we have that 
the condition (t) is equivalent to the condition:
For any function germ $f$ on $(T^* {\mathbb R}^n,0)$,
 there exist a fiber preserving 
function germ $H$ on $(T^* {\mathbb R}^n,0)$ 
and a function germ $g\in B$ such that
$\pi_l(X_f\circ S)=\pi_l( S_*X_g+X_H\circ S)$.
This is equivalent to the condition:
For any function germ $f$ on $(T^* {\mathbb R}^n,0)$,
 there exist a fiber preserving 
function germ $H$ on $(T^* {\mathbb R}^n,0)$ 
and a function germ $g\in B$ such that
$f\circ S-g-H\circ S\in {\mathfrak M}_{T^* {\mathbb R}^n}^{l+2}$.
This is equivalent to the condition:
For any function germ $f$ on $(T^* {\mathbb R}^n,0)$,
 there exist a fiber preserving 
function germ $H$ on $(T^* {\mathbb R}^n,0)$ 
 such that
$f-H\circ S\in B+{\mathfrak M}_{T^* {\mathbb R}^n}^{l+2}$.
This is equivalent to (a').\\
(a')$\Leftrightarrow$(a'').
We need only to prove (a')$\Rightarrow$(a'').
By Margrange preparation theorem, the condition (a') is equivalent that 
 ${\cal E}_{T^* {\mathbb R}^n}/$ $((\pi\circ S)^*{\mathfrak M}_{{\mathbb R}^n}{\cal E}_{T^* {\mathbb R}^n}
+B+
{\mathfrak M}_{T^* {\mathbb R}^n}^{l+2})$ is generated by 
$1,p_1\circ S,\ldots,p_n\circ S$ over ${\mathbb R}$.
This means that 
\[ {\mathfrak M}_{T^* {\mathbb R}^n}^{n+1}\subset 
(\pi\circ S)^*{\mathfrak M}_{{\mathbb R}^n}{\cal E}_{T^* {\mathbb R}^n}
+B+
{\mathfrak M}_{T^* {\mathbb R}^n}^{l+2}.\]
Since ${\mathfrak M}_{T^* {\mathbb R}^n}^{l+2}\subset
{\mathfrak M}_{T^* {\mathbb R}^n}^{n+2}$ ,
it follows  that 
\[ {\mathfrak M}_{T^* {\mathbb R}^n}^{n+1}\subset 
(\pi\circ S)^*{\mathfrak M}_{{\mathbb R}^n}{\cal E}_{T^* {\mathbb R}^n}
+B+
{\mathfrak M}_{T^* {\mathbb R}^n}^{n+2}.\]
Therefore we have that 
\[ {\mathfrak M}_{T^* {\mathbb R}^n}^{n+1}\subset (\pi\circ S)^*{\mathfrak M}_{{\mathbb R}^n}{\cal E}_{T^* {\mathbb R}^n}
+B.\]
This means (a'').\\
(a'')$\Leftrightarrow$(is).
The condition (a'') is equivalent to the condition (a) in Theorem \ref{stable:thleg} and this is equivalent to (is).\\
(t)\&(is)$\Rightarrow$(s).
Since $j^l_0S$ is transversal to $[j^l_0S(0)]$, it follows that 
there exists 
a representative $\tilde{S}\in S(U,T^* {\mathbb R}^n)$ of $S$ and
a neighbourhood  $W_{\tilde{S}}$ of $\tilde{S}$ such that for any $\tilde{S'}\in W_{\tilde{S}}$
 there exists $w\in U$ such that 
$j^l_0\tilde{S}'$ is transversal to  $[j^l_0S(0)]$ at $w$.
Since $j^l_0\tilde{S'_w}(0)\in [j^l_0S(0)]$, it follows that there exists a symplectic embedding germ 
$S''$ on $(T^* {\mathbb R}^n,0)$ such that $\pi\circ i$ and 
$\pi\circ S''|_{\mathbb L}$ is 
Lagrangian equivalent
and $j^l_0S''(0)=j^l_0\tilde{S'_w}(0)$. 
Since $\pi\circ i$ is infinitesimally stable, it follows that $\pi\circ i$ is $l$-determined
by Theorem \ref{n+1det:Lag}.
Therefore we have that $\pi\circ S''|_{\mathbb L}$ is also 
$l$-determined.
It follows that $\pi\circ S''|_{\mathbb L}$ and
 $\pi \circ \tilde{S'_w}|_{\mathbb L}$ is Lagrangian equivalent.
This means that $\pi \circ i$ is stable.
\hfill $\blacksquare$

\vspace{3mm}
Let $S$ be a symplectic diffeomorphism germ on $(T^* {\mathbb R}^n,0),\ 
F(x,y,q)\in {\mathfrak M}(r;k+n)^2$, 
and  $f(x,y)\in {\mathfrak M}(r;k)^2$.
We denote $[S],[F], [f]$ by the equivalence classes of $S,F,f$ under the 
Lagrangian equivalence,
the stably reticular ${\cal P}$-${\cal R}^+$-equivalence, and 
the stably reticular ${\cal R}$-equivalence respectively.
Then the following holds:
\begin{lem}
Let $\pi\circ i_j: ({\mathbb L},0)\rightarrow (T^* {\mathbb R}^n,0)
\rightarrow  ({\mathbb R}^n,0)$ be a stable reticular Lagrangian map with
a generating family $F_j(x,y,q)\in {\mathfrak M}(r;k_j+n)^2$, and 
$S_j$ be an extension of $i_j$ for $j=1,2$.
Then $[S_1]=[S_2]$ if and only if $[F_1]=[F_2]$ and this holds if and only if 
$[F_1|_{q=0}]=[F_2|_{q=0}]$.
\end{lem}
{\it Proof}. 
$[S_1]=[S_2]$ if and only if reticular Lagrangian maps $\pi\circ i_1$ and 
$\pi\circ i_2$ are Lagrangian equivalent by Lemma \ref{SS:lem}, and 
this holds if and only if $[F_1]=[F_2]$
by Theorem \ref{Laggene:th}(3), and 
if this holds then $[F_1|_{q=0}]=[F_2|_{q=0}]$.
Conversely suppose that $[F_1|_{q=0}]=[F_2|_{q=0}]$.
Since $F_1$ and $F_2$ are reticular 
${\cal P}$-${\cal R}^+$-stable unfoldings of $[F_1|_{q=0}]$ and 
$[F_2|_{q=0}]$ respectively and 
$F_1|_{q=0}$ and $F_2|_{q=0}$ are stably reticular 
${\cal R}$-equivalent.
it follows that $[F_1]=[F_2]$.\hfill $\blacksquare$

\vspace{3mm}
\begin{coro}\label{ts:coro}
Let $\pi \circ i: ({\mathbb L},0)\rightarrow (T^* {\mathbb R}^n,0)
\rightarrow  ({\mathbb R}^n,0) $ be a stable reticular Lagrangian map.
Then there exist a neighborhood $U$ of $0$ in  $T^* {\mathbb R}^n$
and $\tilde{S}\in S(U,T^* {\mathbb R}^n)$ with 
$i=\tilde{S_0}|_{\mathbb L}$ \rm{(} \it{that is, $\tilde{S}$ is 
a representative of  a extension of $i$}\rm{)} \it{ such that  
reticular Lagrangian maps  $\pi\circ \tilde{S_w}|_{\mathbb L}$ are
stable for all $w\in U$}
\end{coro}
{\it Proof}.
By Theorem \ref{ts:thlag} (a'), the stability of reticular Lagrangian maps 
are determined by the $(n+1)$-jets of $\pi\circ \tilde{S_w}$ for $w\in U$.
Therefore we have the result by shrinking $U$ if necessary.
\hfill $\blacksquare$

\vspace{3mm}
Let $\pi \circ i: ({\mathbb L},0)\rightarrow (T^* {\mathbb R}^n,0)
\rightarrow  ({\mathbb R}^n,0) $ be a stable reticular Lagrangian map.
We say that $\pi\circ i$ is {\it simple} if 
there exist 
a neighborhood $U$ of $0$ in  $T^* {\mathbb R}^n$
and $\tilde{S}\in S(U,T^* {\mathbb R}^n)$ such that 
$i=\tilde{S_0}|_{\mathbb L}$  and  
$\{ \tilde{S}_w| w\in U\}$ is covered by 
finite orbits $[S_1],\ldots,[S_m]$ for   
symplectic diffeomorphism germs $S_1,\ldots,S_m$ on $(T^* {\mathbb R}^n,0)$.
\begin{lem}\label{simpleLagket:lem}
Let $\pi \circ i: ({\mathbb L},0)\rightarrow (T^* {\mathbb R}^n,0)
\rightarrow  ({\mathbb R}^n,0) $ be a stable reticular Lagrangian map.
Then $\pi \circ i$ is simple if and only if there exist a  neighborhood $U_z$ of 
$z=j^{n+1}_0S(0)$ in $S^{n+1}(n)$ and $z_1,\ldots,z_m\in S^{n+1}(n)$ such 
that 
$U_z\subset [z_1]\cup\cdots\cup [z_m]$.
\end{lem}
{\it Proof}.
Suppose that $\pi\circ i$ is simple.
Then there exists a representative 
$\tilde{S}:U\rightarrow T^* {\mathbb R}^n$ of
 an extension of $i$ and 
symplectic diffeomorphism germs $[S_1],\ldots,[S_m]$ 
on $(T^* {\mathbb R}^n,0)$
such that 
\begin{equation}
\{ \tilde{S}_w| w\in U\}\subset 
[S_1]\cup\cdots\cup [S_m].\label{trans1:eqn}
\end{equation}
Since $\pi\circ i$ is stable, it follows that 
$j^{n+1}_0\tilde{S}$ is transversal to $[z]$ at $0$ by 
Theorem \ref{stable:thlag}.
This means that 
there exists a neighbourhood $U_z$ of $z$ in $S^l(n)$ such that 
$U_z\subset \cup_{w\in U}[j^{n+1}_0\tilde{S}(w)]$.
It follows that  $U_z
\subset 
[j^{n+1}S_1(0)]\cup\cdots\cup [j^{n+1}S_m(0)]$.

Conversely suppose that 
there exist a neighbourhood $U_z$ of $z$
 in $S^{n+1}(n)$
and $z_1,\ldots,z_m
\in S^{n+1}(n)$ such that 
$U_z\subset [z_1]\cup\cdots\cup [z_m]$.
Since the map $j^{n+1}_0\tilde{S}:U\rightarrow S^{n+1}(n)$ is 
continuous,
there exists a neighbourhood $U'$ of $0$ in $U$
such that  $j^{n+1}_0\tilde{S}(w)\in U_z$ for any $w\in {U'}$.
Then we have that 
$\cup_{w\in {U'}}j^{n+1}_0\tilde{S}(w)
\subset [z_1]\cup\cdots\cup [z_m]$.
Choose symplectic diffeomorphism germs $S_1,\ldots,S_m$ on 
$(T^* {\mathbb R}^n,0)$ such that $j^{n+1}S_j(0)=z_j$ 
for $j=1,\ldots,m$.
By Corollary \ref{ts:coro}, we may assume that 
 each $\pi\circ S_j|_{\mathbb L}$ is stable, therefore 
$(n+1)$-determined.
For any $w\in {U'}$ we have that 
there exists $j\in \{1,\ldots,m\}$ such that
$j^{n+1}_0\tilde{S}_w(0)\in [j^{n+1}S_j(0)]$.
It follows that reticular Lagrangian maps 
$\pi \circ \tilde{S_w}|_{\mathbb L}$  and
 $\pi \circ S_j|_{\mathbb L}$
are Lagrangian equivalent.
Therefore $\tilde{S_w}\in [S_j]$.
We have (\ref{trans1:eqn}). \hfill $\blacksquare$
\begin{lem}\label{simpleLag:lem}
A stable reticular Lagrangian map $\pi\circ i$ is simple if and only if 
for a generating family $F(x,y,q)\in {\mathfrak M}(r;k+n)^2$ of 
$\pi\circ i$, 
the function germ $F(x,y,0)\in {\mathfrak M}(r;k)^2$ is ${\cal R}$-simple singularity.
\end{lem}
{\it Proof}.
Let $S$ be an extension of $i$ and $\tilde{S}\in 
 S(U,T^* {\mathbb R}^n)$ be a representative of $S$

Suppose that $F(x,y,0)$ is ${\cal R}$-simple.
The simplicity of  reticular Lagrangian maps is invariant under 
Lagrangian equivalences,
we may assume that the map germ $(Q,P)\mapsto (Q,p\circ S(Q,P))$ is
a diffeomorphism germ on $({\mathbb R}^{2n},0)$.
We consider a symplectic diffeomorphism germ $\tilde{S}_w$ on 
$({\mathbb R}^{2n},0)$ for $w\in U$ near $0$.
Then 
there exists a function germ $H_w(Q,p)\in {\mathfrak M}(2n)^2$
such that the canonical relation $P_w$ associated with $\tilde{S}_w$ has the form:
\[ P_w=\{(Q,-\frac{\partial H_w}{\partial Q}(Q,p),
-\frac{\partial H_w}{\partial p}(Q,p),p)\in 
(T^* {\mathbb R}^n\times T^* {\mathbb R}^n,(0,0))
\}.\]
Then the function germ $F_w(x,y,q)=H'_w(x,y)+\langle y,q\rangle
\in {\mathfrak M}(r;n+n)^2$ is a generating family of $\pi\circ 
\tilde{S}_w|_{\mathbb L}$,
where $H'_w\in {\mathfrak M}(r;n)^2$ is defined by 
$H'_w(x,y)=H_w(x,0,y)$.
Since $F_0$ is a generating family of $\pi\circ i$, 
we have that $F_0(x,y,0)(=H'_0(x,y))$ is stably ${\cal R}$-equivalent to $F(x,y,0)$.
Therefore we have that $H'_0$ is 
${\cal R}$-simple.
Then there exists $f_1,\ldots,f_m\in {\mathfrak M}(r;n)$ and 
a neighbourhood $V$ of $j^{n+2}H'_0(0)$ in $J^{n+2}(r+n,1)$ 
such that $V\subset [j^{n+2}f_1(0)]\cup\cdots\cup [j^{n+2}f_m(0)]$.

Since  the $(n+2)$-jet of $H'_S$ is determined by the $(n+1)$-jet of $S$,
there exists a neighbourhood $U'$ of $0$ in $U$
such that the map germ
\[ U'\rightarrow J^{n+2}(r+n,1),\ w \mapsto j^{n+2}H'_w(0)\]
is well defined and continuous.

Let $U''$ be the inverse image of $V$ by the above map.
Then for any $w\in {U''}$ 
the reticular Lagrangian map $\pi\circ \tilde{S}_w|_{\mathbb L}$ has 
a generating family which is reticular ${\cal P}$-${\cal R}^+$-equivalent to
$f_j(x,y)+\langle y,q\rangle\in {\mathfrak M}(r;n+n)^2$ for some $j$
because $\pi\circ \tilde{S}_w|_{\mathbb L}$ is stable by Corollary 
\ref{ts:coro} and hence $(n+1)$-determined.
It follow that $\{ \tilde{S}_w |w\in {U''}\}\subset [S_1]\cup\cdots\cup
[S_m]$,
where $S_j$ is an extension of a reticular Lagrangian embedding which
defines a reticular Lagrangian map with the generating family 
$f_j(x,y)+\langle y,q\rangle\in {\mathfrak M}(r;n+n)^2$.
This means that $\pi\circ i$ is simple.

Conversely suppose that  $\pi\circ i$ is simple.
Let $S_0$ be an extension of $i$.
we may assume that 
there exists a function germ $H_0(Q,p)\in {\mathfrak M}(2n)^2$
such that the canonical relation $P_{S_0}$ associated with $S_0$ 
has the form:
\[ P_{S_0}=\{(Q,-\frac{\partial H_0}{\partial Q}(Q,p),
-\frac{\partial H_0}{\partial p}(Q,p),p)\in 
(T^* {\mathbb R}^n\times T^* {\mathbb R}^n,(0,0))
\}.\]
This means that for a function germ $H$ around $H_0$ 
there exists a symplectic diffeomorphism germ $S_H$ on 
$(T^* {\mathbb R}^n,0)$ with the canonical relation $P_{S_H}$
which has the form same as $P_{S_0}$.
Then the $(n+1)$-jet of $S_H$ is determined by the $(n+2)$-jet 
of $H$.
The function germ 
$H'_0(x,0,y)+\langle y,q\rangle\in {\mathfrak M}(r;n+n)^2$ is 
a generating family of $\pi\circ i$,
where $H'_0\in {\mathfrak M}(r;n)^2$ is defined by 
$H'_0(x,y)=H_0(x,0,y)$.
Then there exists a quadratic form $T(Q_{r+1},\ldots,Q_n,
p_1,\ldots,p_n)$ such that 
the function germ $H''_0(Q,p)=H'_0(Q_1,\ldots,Q_r,p)+
T(Q_{r+1},\ldots,Q_n,p)$ also defined the 
symplectic diffeomorphism $S_{H''_0}$.
Then we have that $\pi\circ S_{H''_0}|_{\mathbb L}=\pi\circ i$.
Therefore we may assume that $H''_0=H_0$

We have that $H'_0$ is stably 
reticular ${\cal R}$-equivalent to $F(x,y,0)$.
Therefore we  need only to prove that $H'_0$
is a simple singularity.

There exists a neighbourhood $V$ of $j^{n+2}H'_0(0)$ such that 
the map 
\[ V\rightarrow S^{n+1}(n),\ j^{n+2}H'(0)\mapsto j^{n+1}_0S_H(0),\]
where $H(Q,p)=H'(Q_1,\ldots,Q_r,p)+
T(Q_{r+1},\ldots,Q_n,p)$ 
is well defined and continuous.

Since $\pi\circ i$ is simple,
there exist finite symplectic diffeomorphism germs 
$S_1,\ldots,S_m$ on  $(T^* {\mathbb R}^n,0)$ and
a neighbourhood $U_z$ of $z=j^{n+1}S_0(0)$ in $S^{n+1}(n)$ such that
$U_z\subset [S_1]\cup\cdots\cup [S_m]$
Let $V'$ be the inverse image of $U_z$ by the above map.
Then we have that $V'\subset [z_1]\cup\cdots\cup [z_m]$,
where $z_j=j^{n+2}H_j(0)$ for a function germ $H'_j\in 
{\cal E}(r;n)$ such that $j^{n+1}S_j(0)=j^{n+1}S_{H_j}(0)$.
This means that $H'_0$ is simple.

\begin{tth}\label{genericclassLag:th}
Let $r=0,n\leq 5$ or $r=1,n\leq 3$. 
Let  $U$ be a neighborhood of $0$ in $T^* {\mathbb R}^n$.
Then there exists a residual set $O\subset S(U,T^* {\mathbb R}^n)$
such that for any $\tilde{S}\in O$ and $w=(0,\ldots,p^0_{r+1},\ldots,p^0_n)\in U$,
the reticular Lagrangian map $\pi\circ \tilde{S}_w|_{\mathbb L}$ 
is stable.
\end{tth}

In the case $r=0,\ n\leq 5$. 
A reticular Legendrian map $\pi\circ \tilde{S}_w|_{\mathbb L}$ 
for any $\tilde{S}\in O$ and $w\in U$
has a generating family $F$ which is a 
reticular ${\cal P}$-${\cal R}^+$-stable unfolding of one of 
$A_2,A^\pm_3,A_4,A^\pm_5,A_6,
D^\pm_4,D^\pm_4,D^\pm_6,E^\pm_6$, 
that is $F$ is stably reticular ${\cal P}$-${\cal R}^+$-equivalent to 
one of the following list:\\
$A_2:F(y_1,q_1)=y_1^3+q_1y_1$,\\
$A^\pm_3:F(y_1,q_1,q_2)=\pm y_1^4+q_1y_1^2+q_2y_1$,\\
$A_4:F(y_1,q_1,q_2,q_3)= y_1^5+q_1y_1^3+q_2y_1^2+q_3y_1$,\\
$A^\pm_5:F(y_1,q_1,q_2,q_3,q_4)= 
\pm y_1^6+q_1y_1^4+q_2y_1^3+q_3y_1^2+q_4y_1$,\\
$A_6:F(y_1,q_1,q_2,q_3,q_4,q_5)= 
y_1^7+q_1y_1^5+q_2y_1^4+q_3y_1^3+q_4y_1^2+q_5y_1$,\\
$D^\pm_4:F(y_1,y_2,q_1,q_2,q_3)=
y_1^2y_2\pm y_2^3+q_1y_2^2+q_2y_2+q_3y_1$,\\
$D^\pm_5:F(y_1,y_2,q_1,q_2,q_3,q_4)=
y_1^2y_2\pm y_2^4+q_1y_2^3+q_2y_2^2+q_3y_2+q_4y_1$,\\
$D^\pm_6:F(y_1,y_2,q_1,q_2,q_3,q_4,q_5)=
y_1^2y_2\pm  y_2^5+q_1y_2^4+q_2y_2^3+q_3y_2^2+q_4y_2
+q_5y_1$,\\
$E^\pm_6:F(y_1,y_2,q_1,q_2,q_3,q_4,q_5)=
y_1^3\pm  y_2^4+q_1y_1y_2^2+q_2y_1y_2+q_3y_2^2+q_4y_1
+q_5y_2$.
\vspace{3mm}\\
In the case $r=1,\ n\leq 3$.   
A reticular Legendrian map $\pi\circ \tilde{S}_w|_{\mathbb L}$
for any $\tilde{S}\in O$ and $w\in U$
has a generating family which is a 
${\cal P}$-${\cal R}^+$-stable unfolding of one of 
$B^\pm_2,B^\pm_3,B^\pm_4,C^\pm_3,C^\pm_4,F^\pm_4$,
that is $F$ is stably reticular ${\cal P}$-${\cal R}^+$-equivalent to 
one of the following list:\\
$B^\pm_2:F(x,q_1)=\pm x^2+q_1x$,\\
$B^\pm_3:F(x,q_1,q_2)=\pm x^3+q_1x^2+q_2x$,\\
$B^\pm_4:F(x,q_1,q_2,q_3)=\pm x^4+q_1x^3+q_2x^2+q_1x$,\\
$C^\pm_3:F(x,y,q_1,q_2)=\pm xy+y^3+q_1y^2+q_2y$,\\
$C^\pm_4:F(x,y,q_1,q_2,q_3)=\pm xy+y^4+q_1y^3+q_2y^2+q_3y$,\\
$F^\pm_4:F(x,y,q_1,q_2,q_3)=\pm x^2+y^3+q_1xy+q_2x+q_3y$.
\vspace{3mm}\\
{\it Proof}. We need only to prove the case $r=1,\ n\leq 3$.
Let $F_X(x,y,q)\in {\mathfrak M}(r;k+n)^2$ 
be a reticular ${\cal P}$-${\cal R}^+$-stable unfolding of singularity 
$X\in {\mathfrak M}(r;k)^2$ for 
\[ X=B^\pm_2,B^\pm_3,B^\pm_4,
C^\pm_3,C^\pm_4,F^\pm_4.\]
Then other unfoldings are not stable since other singularities have 
reticular ${\cal R}^+$-codimension $>3$.
We choose stable reticular Lagrangian maps $\pi\circ i_X:
({\mathbb L},0)\rightarrow (T^* {\mathbb R}^n,0)
\rightarrow  ({\mathbb R}^n,0)$ with generating family $F_X$
 and $S_X$ be and extension of $i_X$ for above list.
We set $L_U=\{ (q,p)\in U |q=p_1=0 \}$ and define that
\[ O'=\{ \tilde{S}\in S(U,T^* {\mathbb R}^n)\ |\ 
j^{n+1}_0\tilde{S} \mbox{ is transversal to }[j^{n+1}S_X(0)]
\mbox{ on } L_U \mbox{ for all } X \}.\]
Then $O'$ is a residual set. We set 
\[ Y=\{ j^{n+1}S(0)\in S^{n+1}(n)\ |\ 
\mbox{the codimension of }[j^{n+1}S(0)]>2n\}.\]
Then $Y$ is an algebraic set in $S^{n+1}(n)$ by Theorem \ref{ts:thlag} (a').
Therefore we can define that 
\[ O''=\{ \tilde{S}\in S(U,T^* {\mathbb R}^n)\ |\ 
j^{n+1}_0\tilde{S} \mbox{ is transversal to } Y \}.\]
Then $Y$ has codimension $>2n$ because all element in $Y$ is 
adjacent to one of the above list which are simple.
Then we have that
\[ O''=\{ \tilde{S}\in S(U,T^* {\mathbb R}^n)\ |\ 
j^{n+1}_0\tilde{S}(U)\cap Y=\emptyset \}.\]
We define $O=O'\cap O''$. Then $O$ has the required condition.

In the case $r=0,n\leq 5$. 
Set $X=A_2,A^\pm_3,A_4,A^\pm_5,A_6,
D^\pm_4,D^\pm_4,D^\pm_6,E^\pm_6$ and 
$Y=\{ j^{n+1}S(0)\in S^{n+1}(n)\ |\ 
\mbox{the codimension of }[j^{n+1}S(0)]>2n\}$.
Then we have that the codimension of $Y$ in $S^{n+1}(n)$ is 
higher than $2n$ and the assertion is proved by the parallel method of 
the above case.
\hfill $\blacksquare$
\begin{figure}[htbp]
\[
\begin{array}{cc}
\includegraphics*[width=6cm,height=5cm]{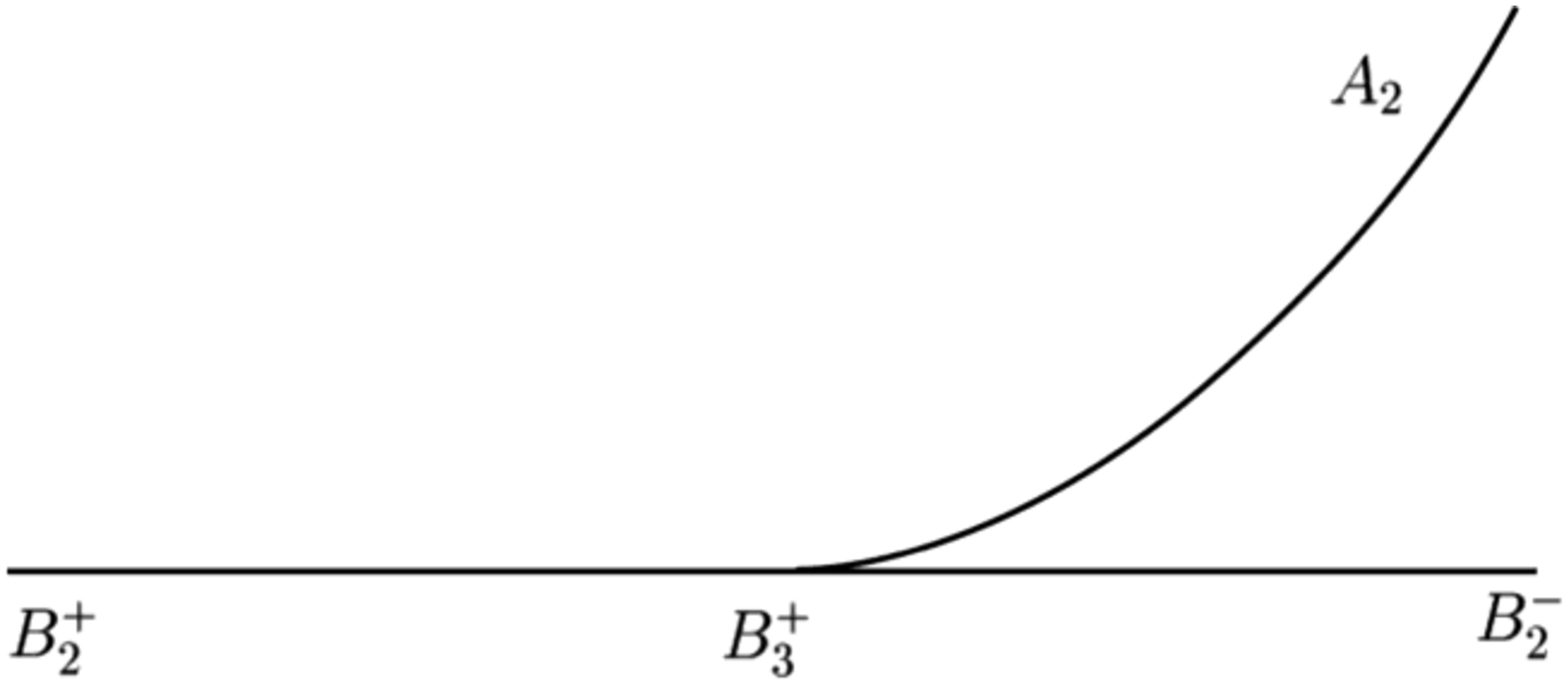}
&
\includegraphics*[width=6cm,height=5cm]{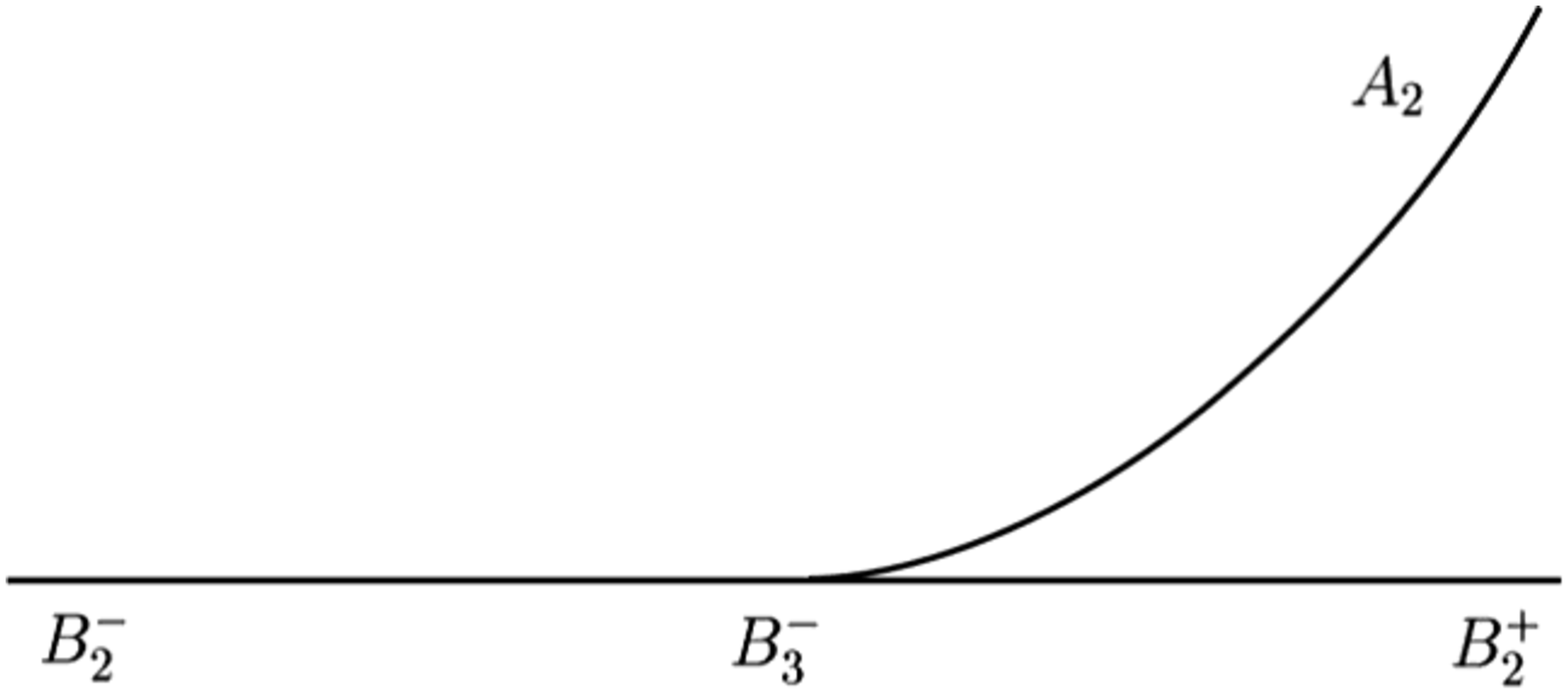}
\end{array}\]
\caption{the caustic $B^\pm_3$}
\[
\begin{array}{cc}
\includegraphics*[width=6cm,height=5cm]{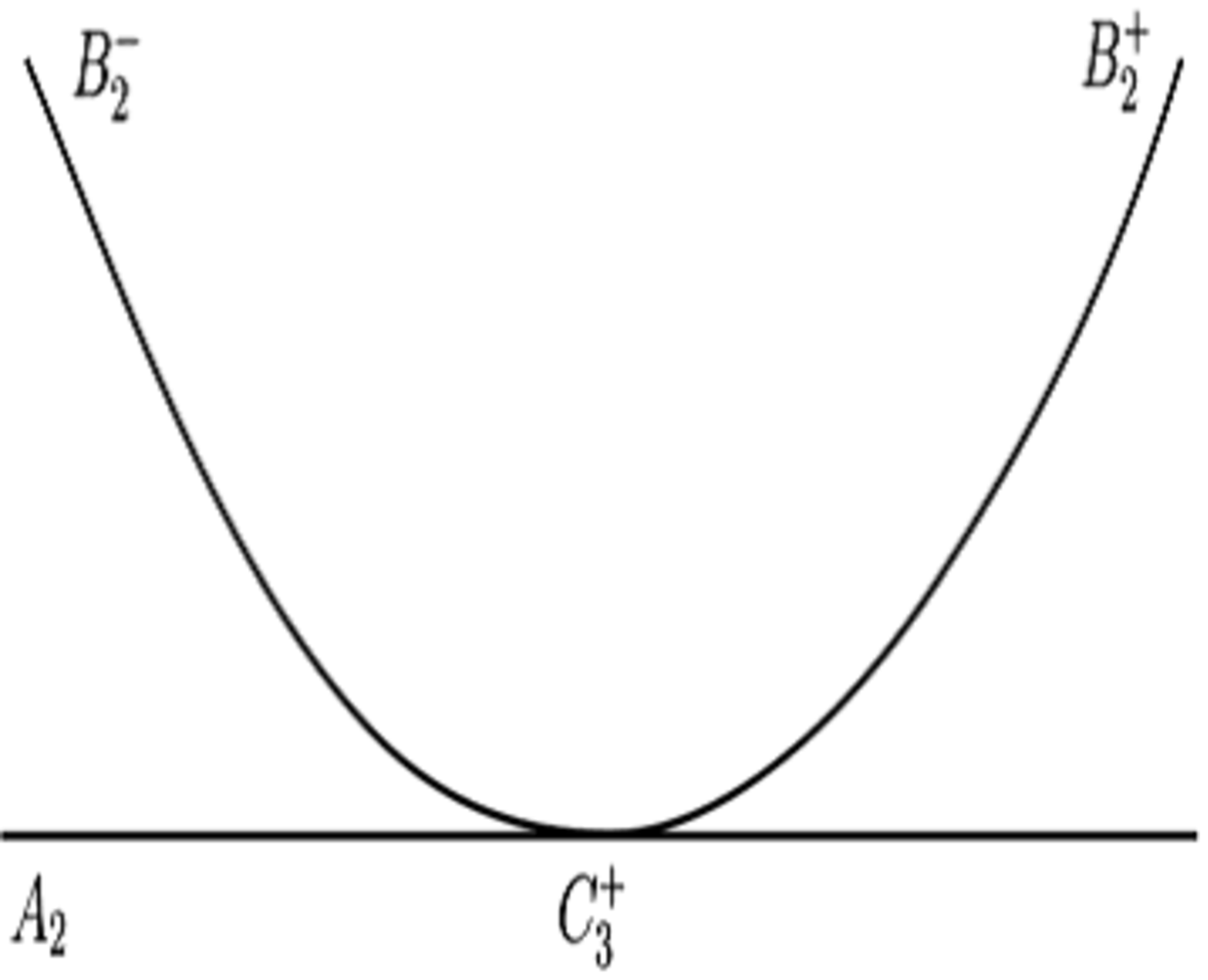}
&
\includegraphics*[width=6cm,height=5cm]{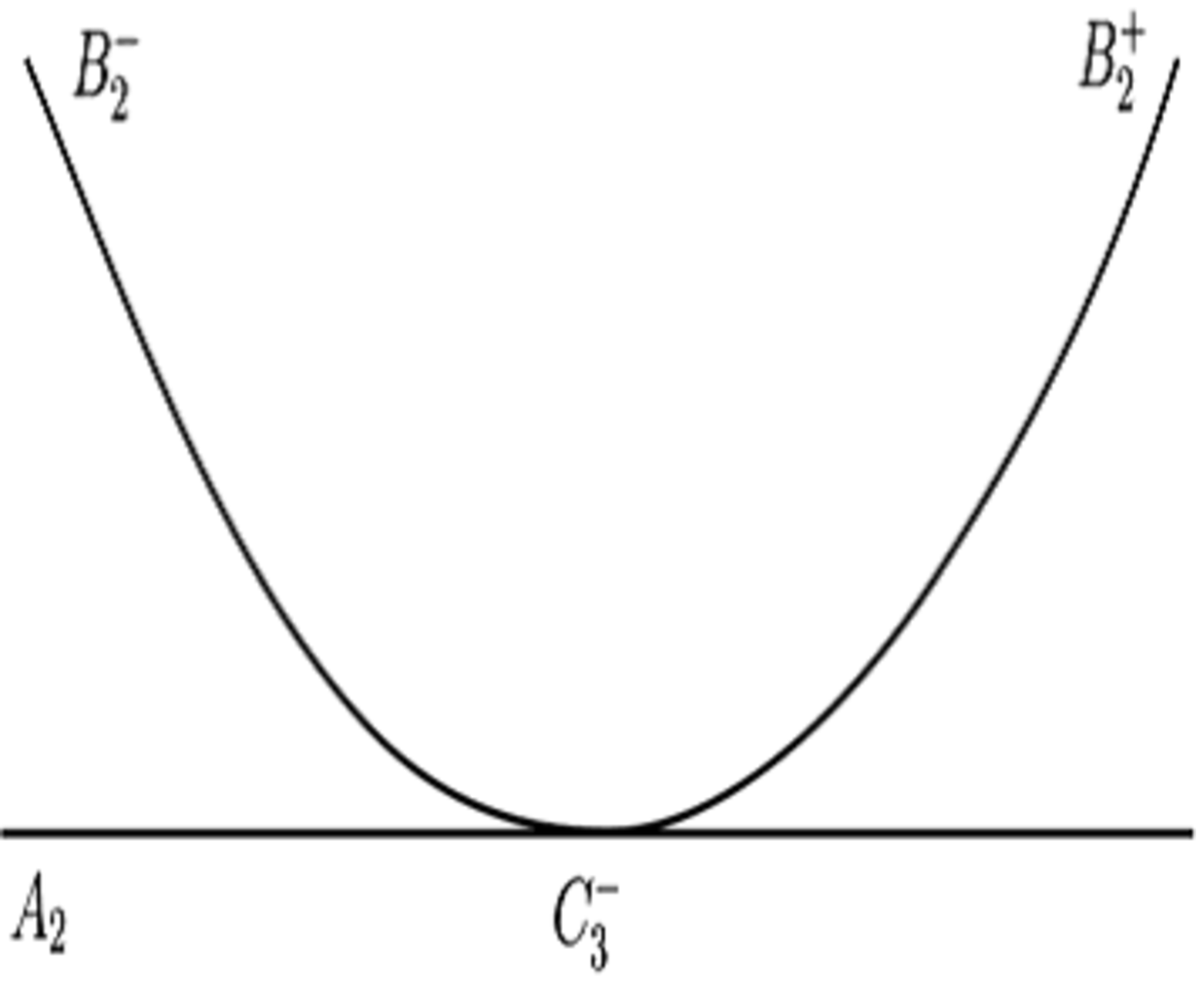}
\end{array}\]
\caption{the caustics $C^\pm_3$}
\end{figure}%

\begin{figure}[htbp]
\[
\begin{array}{cc}
\includegraphics*[width=7cm,height=7cm]{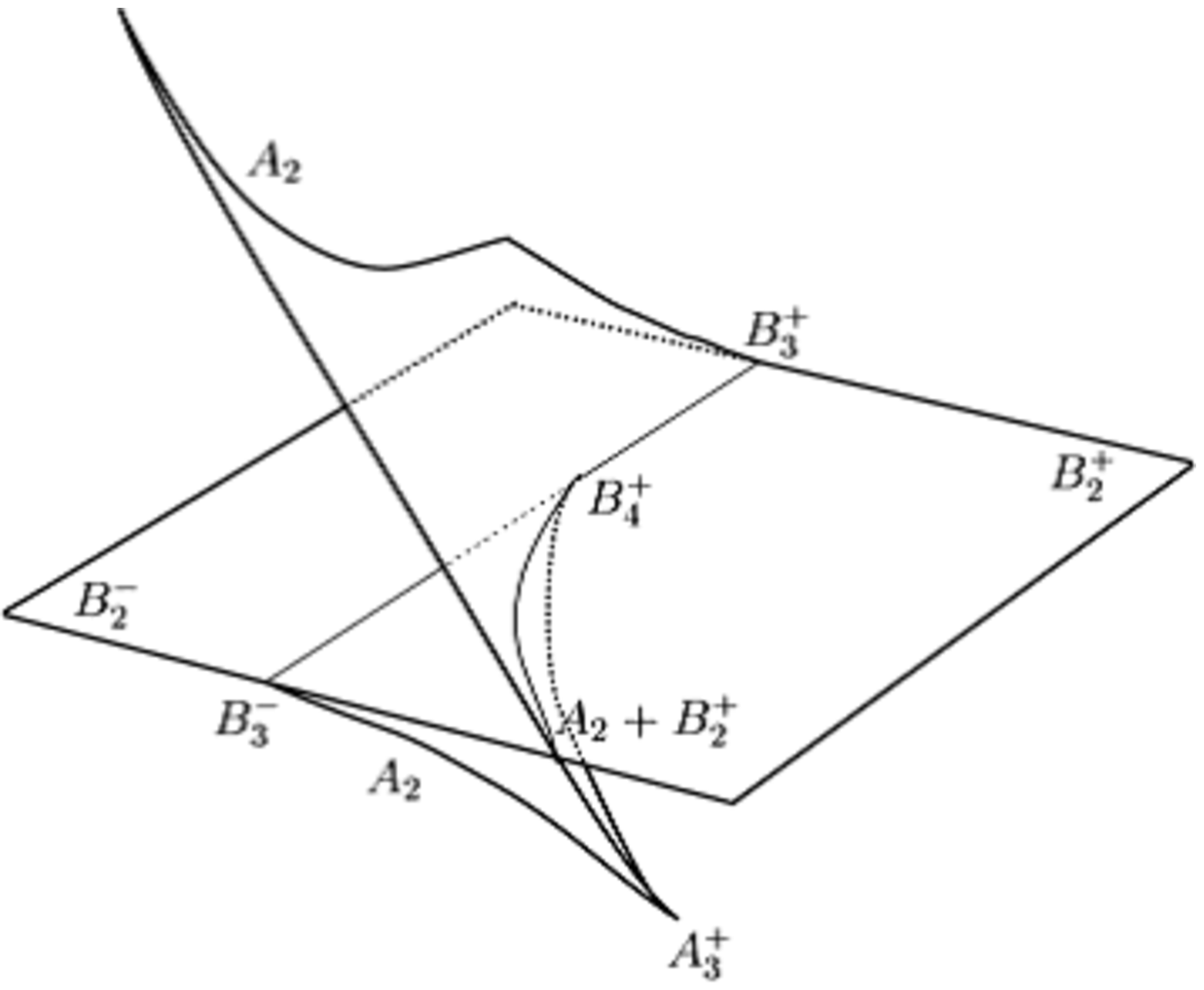}
&
\includegraphics*[width=7cm,height=7cm]{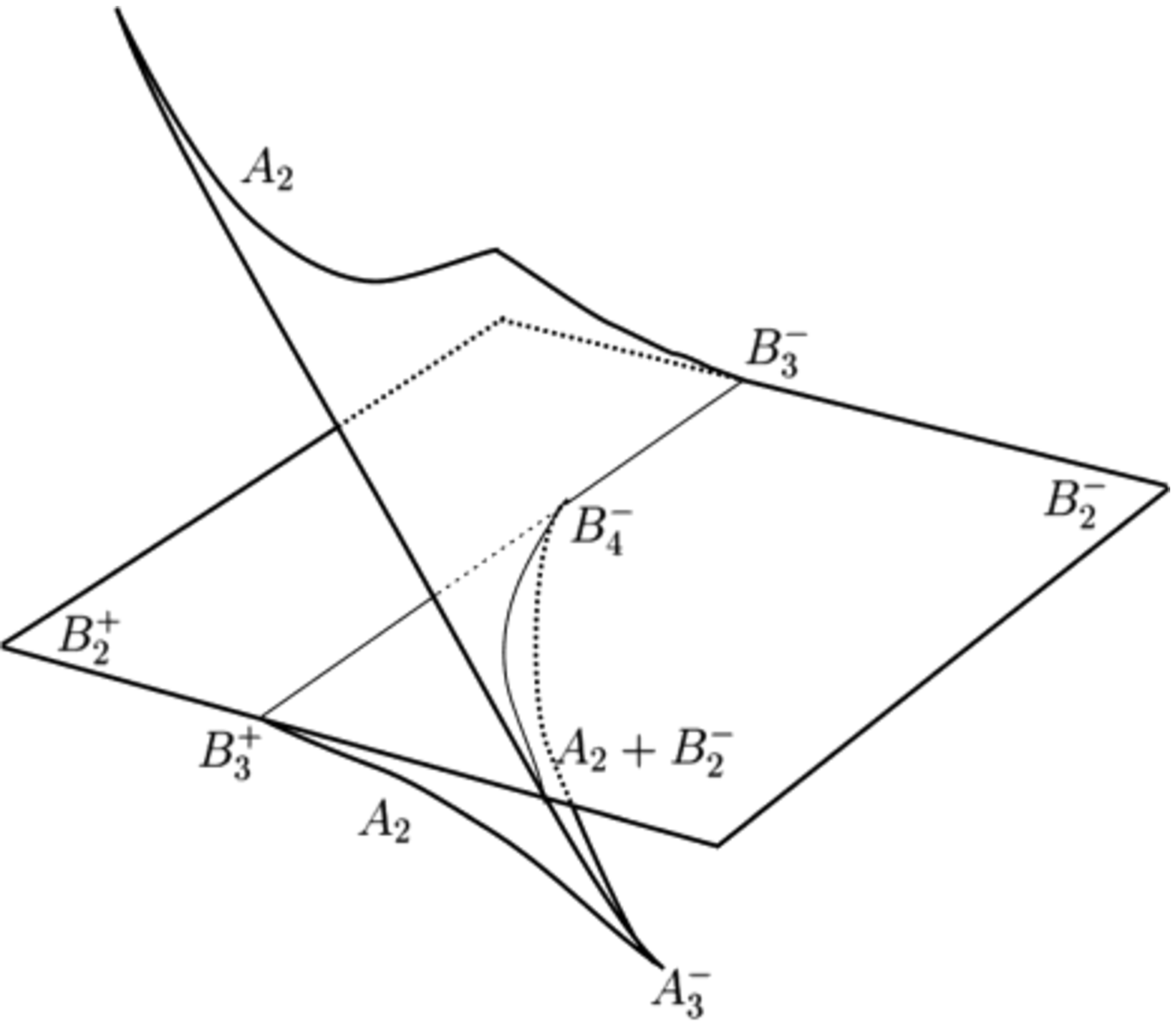}
\end{array}\]
\caption{the caustics $B^\pm_4$}
\[
\begin{array}{cc}
\includegraphics*[width=7cm,height=7cm]{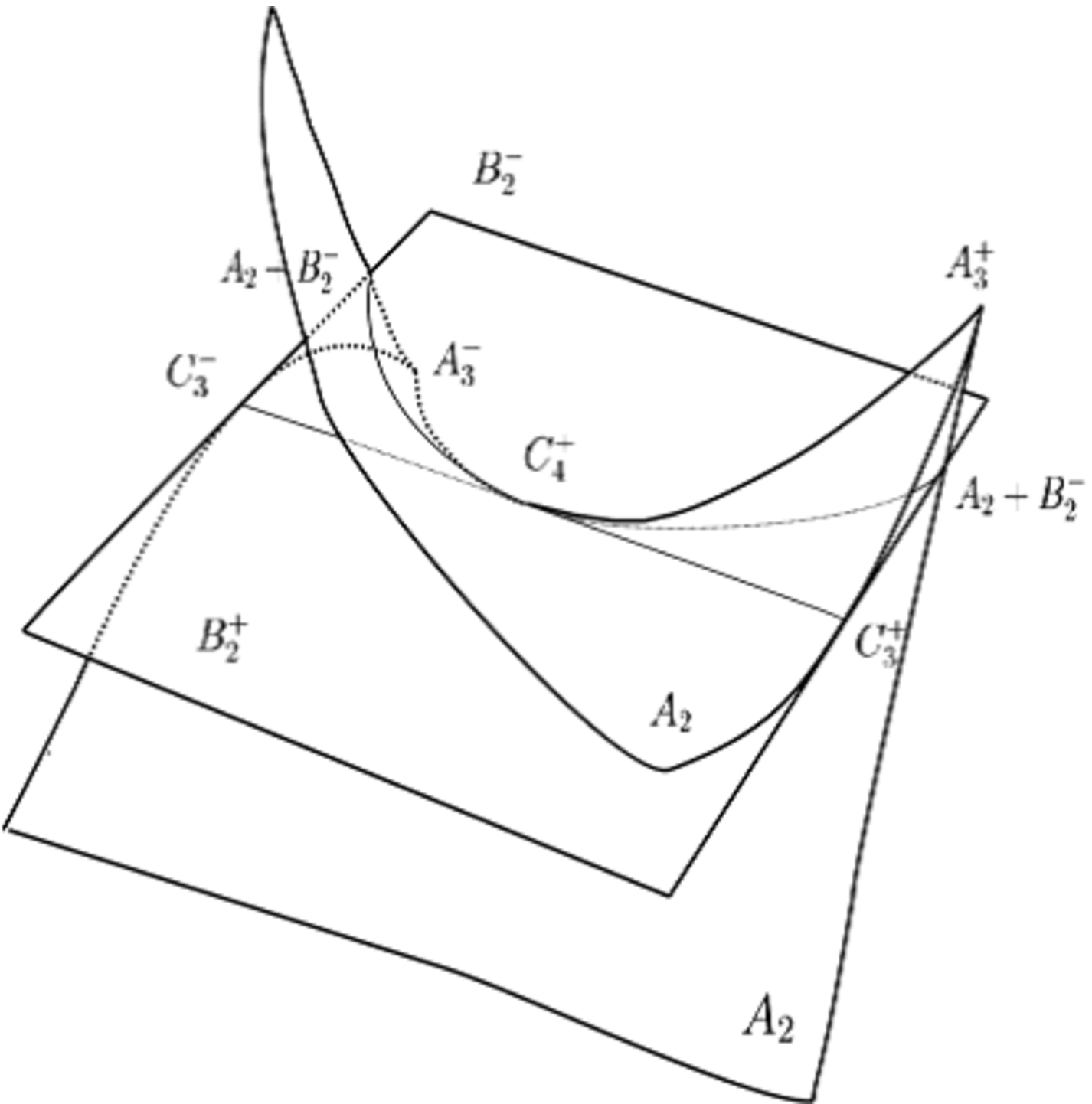}
&
\includegraphics*[width=7cm,height=7cm]{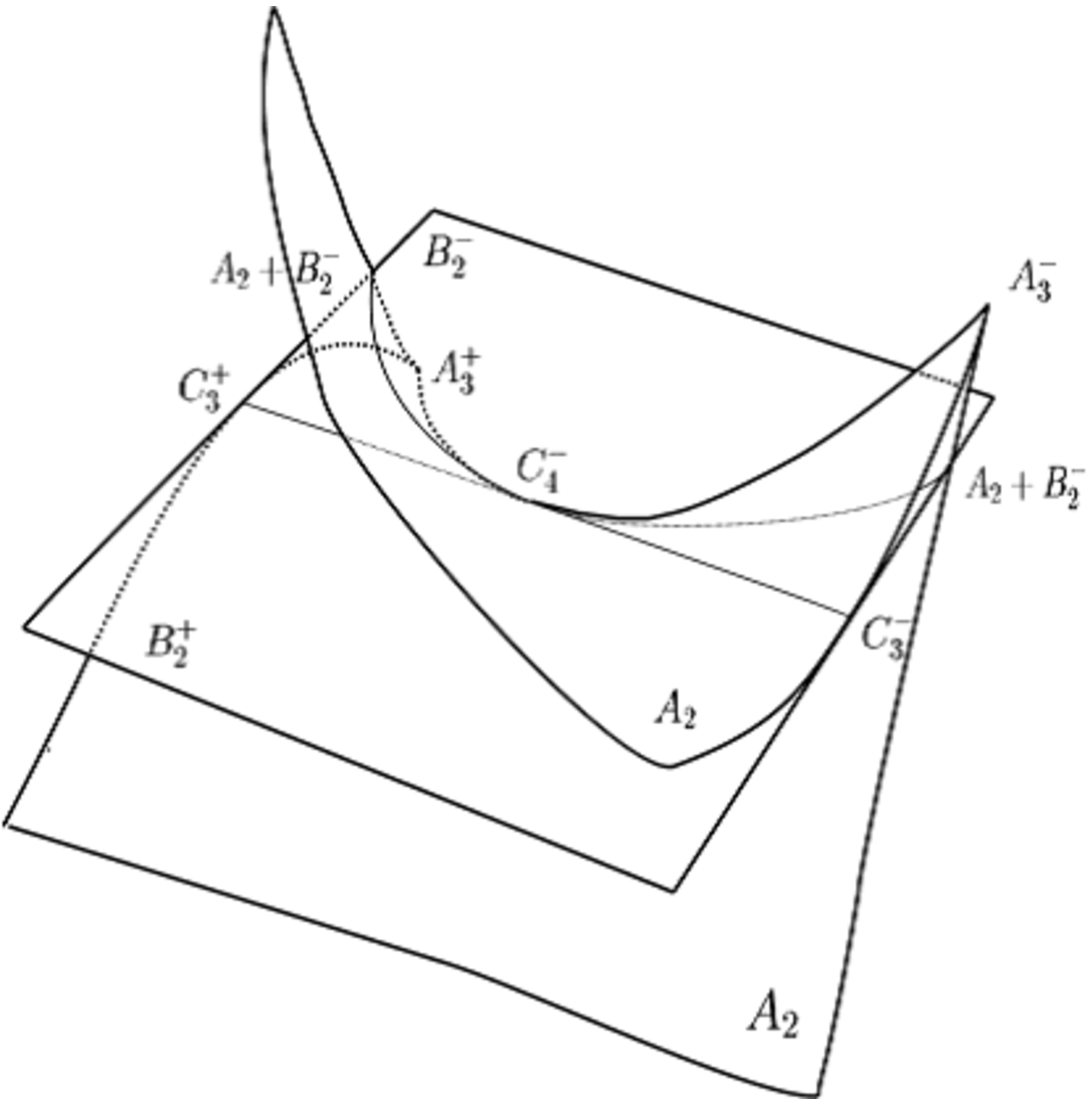}
\end{array}\]
\caption{the caustics $C^\pm_4$}
\end{figure}%

\begin{figure}[htbp]
\[
\begin{array}{cc}
\includegraphics*[width=7cm,height=7cm]{F+_4.eps}
&
\includegraphics*[width=7cm,height=7cm]{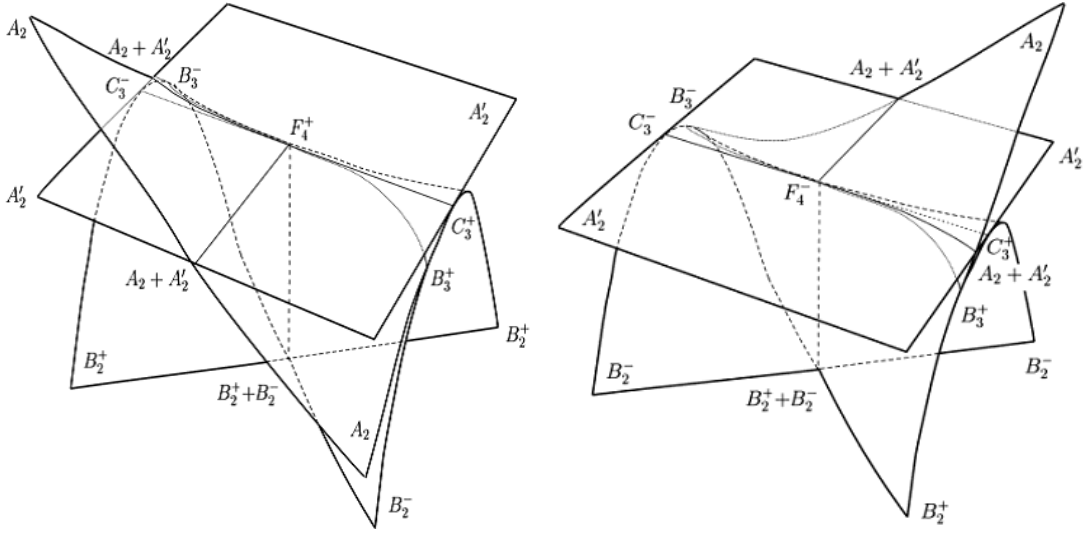}
\end{array}\]
\caption{the caustics $F^\pm_4$}
\end{figure}%
\newpage
\section{Genericity of reticular Legendrian maps}\label{GenLeg}
\quad
Let  $J^l(2n+1,2n+1)$ be the set of $l$-jets of map germs  from $(J^1({\mathbb R}^n,{\mathbb R}),0)$ 
to $(J^1({\mathbb R}^n,{\mathbb R}),0)$ and 
$C^l(n)$ be the Lie group in $J^l(2n+1,2n+1)$ 
consists of $l$-jets of contact diffeomorphism germs 
on $(J^1({\mathbb R}^n,{\mathbb R}),0)$.
We  consider the Lie subgroup $rLe^l(n)$ of $L^l(2n+1)\times L^l(2n+1)$ consists of 
 $l$-jets of reticular diffeomorphisms on the  source space and $l$-jets of Legendrian equivalences of $\tilde{\pi}$: 
\begin{eqnarray*}
rLe^l(n)=\{ (j^l\phi(0),j^l\Theta(0))\in L^l(2n+1)\times L^l(2n+1)\ | \ \phi
 \mbox{ is a reticular}\hspace{3cm}\\ \mbox{ diffeomorphism on }
(J^1({\mathbb R}^n,{\mathbb R}),0), \Theta \mbox{ is a Legendrian equivalence of } 
\tilde{\pi}\}. 
\end{eqnarray*}
The group $rLe^l(n)$ acts on $J^l(2n+1,2n+1)$ and $C^l(n)$ is invariant 
under this action.
Let $C$ be  a contact diffeomorphism germ on $(J^1({\mathbb R}^n,{\mathbb R}),0)$ and set 
$z=j^lC(0)$.
We denote  the orbit $rLe^l(n)\cdot z$ by $[z]$.
Then 

\[ [z]=\{ j^lC'(0)\in C^l(n)\ | \ \tilde{\pi} \circ i \mbox{ and } \tilde{\pi}\circ C'|_{\tilde{\mathbb L}} 
\mbox{ are Legendrian equivalent} \}. \]

In this section we denote by $X_f$ the Contact Hamiltonian vector field on 
$(J^1({\mathbb R}^n,{\mathbb R}),0)$ for a function germ $f$ on 
$(J^1({\mathbb R}^n,{\mathbb R}),0)$.
That is
\[ X_f=\sum_{j=1}^n(\frac{\partial f}{\partial q_j}
+p_j\frac{\partial f}{\partial z})\frac{\partial }{\partial p_j}
-\sum_{j=1}^n\frac{\partial f}{\partial p_j}\frac{\partial }{\partial q_j}
+(f-\sum_{j=1}^np_j\frac{\partial f}{\partial p_j})\frac{\partial }{\partial z}.\]
We denote by $VI_C$  the vector space consists of infinitesimal contact transformation germs of $C$
and  denote by $VI^0_C$ the subspace of $VI_C$ consists of germs which
 vanish on $0$.
We denote by $VL_{J^1({\mathbb R}^n,{\mathbb R})}$ by the vector space consists of infinitesimal Legendrian 
equivalences on $(J^1({\mathbb R}^n,{\mathbb R}),0)$ and denote by $VL^0_{J^1({\mathbb R}^n,{\mathbb R})}$
by the subspace of $VL_{J^1({\mathbb R}^n,{\mathbb R})}$ consists of germs which vanish at $0$.
We denote by $V^0_{\tilde{\mathbb L}}$ the vector space consists of infinitesimal reticular diffeomorphisms 
on $(J^1({\mathbb R}^n,{\mathbb R}),0)$ which vanishes at $0$:
\[ V^0_{\tilde{\mathbb L}}=\{ \xi \in X(J^1({\mathbb R}^n,{\mathbb R}),0)\ | \ 
\xi \mbox{ is tangent to } \tilde{L}^0_\sigma \mbox { for all } \sigma\subset I_r,\ \xi(0)=0\}. \]
\begin{lem}\label{infsta:Leglem}
{\rm (1)} A vector field germ $v$ on $(J^1({\mathbb R}^n,{\mathbb R}),0)$ along $C$ is 
belongs to $VI_C$ if and only if there exists a function germ $f$ on $(J^1({\mathbb R}^n,{\mathbb R}),0)$
such that $v=X_f\circ C$,\\
{\rm (2)} A vector field germ $\eta$ on $(J^1({\mathbb R}^n,{\mathbb R}),0)$ is
belongs to $VL_{J^1({\mathbb R}^n,{\mathbb R})}$ if and only if there exists a fiber preserving function germ
$H$ on $(J^1({\mathbb R}^n,{\mathbb R}),0)$ such that $\eta=X_H$.\\
{\rm (3)}  A vector field germ $\xi$ on $(J^1({\mathbb R}^n,{\mathbb R}),0)$ is
belongs to $V^0_{\tilde{\mathbb L}}$ if and only if there exists a function germ $g\in B'_0$
such that $\xi=X_g$, where $B'_0=\langle q_1p_1,\ldots,q_rp_r\rangle_{{\cal E}_{J^1({\mathbb R}^n,{\mathbb R})}}+
{\mathfrak M}_{J^1({\mathbb R}^n,{\mathbb R})}\langle q_{r+1},\ldots,
$ $q_n,z\rangle$ is a submodule of ${\cal E}_{J^1({\mathbb R}^n,{\mathbb R})}$
\end{lem}
By this lemma we have that:
\[ VI^0_C=\{ v:(J^1({\mathbb R}^n,{\mathbb R}),0)\rightarrow (T(J^1({\mathbb R}^n,{\mathbb R})),0)\ | \
v=X_f\circ C \mbox{ for some }f\in {\mathfrak M}^2_{J^1({\mathbb R}^n,{\mathbb R})}\},\]
\begin{eqnarray*}
VL^0_{J^1({\mathbb R}^n,{\mathbb R})}=\{ \eta\in X(J^1({\mathbb R}^n,{\mathbb R}),0)\ | \
\eta=X_H \hspace{7cm}\\
\mbox{ for some fiver preserving function germ } H\in {\mathfrak M}^2_{J^1({\mathbb R}^n,{\mathbb R})}\},
\end{eqnarray*}
\[ V^0_{\tilde{\mathbb L}}=\{ \xi \in X(J^1({\mathbb R}^n,{\mathbb R}),0)\ | \ 
\xi=X_g  \mbox{ for some }g\in B'_0\}. \]

\vspace{2mm}
We define the homomorphism $tC:
VI^0_{\tilde{\mathbb L}}\rightarrow VI^0_C$ by 
$tC(v)=C_*v$ and define the homomorphism 
$wC:VL^0_{J^1({\mathbb R}^n,{\mathbb R})}\rightarrow 
VI^0_C$ by 
$wC(\eta)=\eta\circ C$.

\vspace{2mm}
We denote $VI^{l}_C$ the subspace of $VI_C$ consists of infinitesimal
contact transformation germs of $C$ whose $l$-jets are $0$:
\[  VI^l_C=\{ v\in VI_C\ | \ j^lv(0)=0\}. \]

\vspace{3mm}
We have that $j^l_0C$ is transversal to $[z]$ if and only if
\[ 
tC( V^0_{\tilde{\mathbb L}})+wC(VL_{J^1({\mathbb R}^n,{\mathbb R})}))+VI^{l+1}_C=VI_C. \]

\vspace{5mm}

Let $N,M$ be $(2n+1)$-dimensional contact manifolds.
We denote $C(N,M)$ the space of contact embeddings from
$N$ to $M$ with the topology induced from the Whitney $C^\infty$-topology
of $C^\infty(N,M)$
and define that
 \[ J^l_C(N,M)=\{ j^lC(u_0)\in J^l(N,M)| C:(N,u_0)\rightarrow M \mbox{ 
is a contact embedding germ},\ u_0\in N\}.\]
\begin{prop}
$C(N,M)$ is a Baire space.
\end{prop}
\begin{tth}[Contact transversality theorem]\label{contacttrans:tth}
Let $N,M$ be $(2n+1)$-dimensional contact manifolds.
 Let 
$Q_j, j=1,2,\ldots$ be submanifolds of $J^l_C(N,M)$.
Then the set
\[ T=\{ C\in C(N,M)\ | \ j^lC \mbox{ is transversal to } Q_j 
\mbox{ for all } j\in {\mathbb N}\} \]
is residual set in $C(N,M)$. In particular $T$ is dense.
\end{tth}
\begin{tth}\label{transv:leg}
Let $U$ be a neighborhood $0$ in $J^1({\mathbb R}^n,{\mathbb R})$, 
$Q_1,Q_2,\ldots$ are submanifolds of $C^l(n)$. 
We define $L_U=\{ (q,z,p)\in U |q=z=p_1=\cdots =p_r=0 \}$
Then the set
\begin{eqnarray*}
T=\{ \tilde{C}\in C(U,J^1({\mathbb R}^n,{\mathbb R}))| j^l_0 \tilde{C}
 \mbox{ is transversal to } Q_j \mbox{ on } L_U
\mbox{ for all } j \}
\end{eqnarray*}
is a residual set in $C(U,J^1({\mathbb R}^n,{\mathbb R}))$.
\end{tth}
We have the following theorem which is proved by a parallel method of Theorem \ref{ts:thlag}
\begin{tth}
Let $\tilde{\pi}\circ i:(\tilde{\mathbb L},0) \rightarrow (J^1({\mathbb R}^n,{\mathbb R}),0) 
\rightarrow({\mathbb R}^n\times{\mathbb R},0)$ be a reticular Legendrian map. 
Let $C$ be an extension of $i$ and $l\geq n+2$.
Let $B'=\langle q_1p_1,\ldots,q_rp_r,q_{r+1},\ldots,q_n,z\rangle_{{\cal E}_{J^1({\mathbb R}^n,{\mathbb R})}}$
be a submodule of ${\cal E}_{J^1({\mathbb R}^n,{\mathbb R})}$.
Then the followings are equivalent:\\
{\rm (s)} $\tilde{\pi}\circ i$ is stable.\\
{\rm (t)} $j^l_0C$ is transversal to $[j^l_0C(0)]$.\\
{\rm (a')} ${\cal E}_{J^1({\mathbb R}^n,{\mathbb R})}/
(B'+{\mathfrak M}_{J^1({\mathbb R}^n,{\mathbb R})}^{l+2})$ is generated by 
$1,p_1\circ C,\ldots,p_n\circ C$ as ${\cal E}_{q,z}$-module via $\tilde{\pi}\circ C$.\\
{\rm (a'')} ${\cal E}_{J^1({\mathbb R}^n,{\mathbb R})}/B'$ is generated by 
$1,p_1\circ C,\ldots,p_n\circ C$ as ${\cal E}_{q,z}$-module via $\tilde{\pi}\circ C$.\\
{\rm (is)} $\tilde{\pi}\circ i$ is infinitesimally stable. 
\end{tth}

\vspace{3mm}
For $\tilde{C}\in C(U,J^1({\mathbb R}^n,{\mathbb R}))$ we define
 the continuous map 
$j^l_0\tilde{C}:U\rightarrow C^l(n)$ by $w$ to the $l$-jet of 
$\tilde{C}_w$.

\vspace{3mm}
Let $\tilde{\pi} \circ i: (\tilde{\mathbb L},0)\rightarrow 
(J^1({\mathbb R}^n,{\mathbb R}),0)
\rightarrow  ({\mathbb R}^n\times{\mathbb R},0)$ 
be a stable reticular Legendrian map.
We say that $\tilde{\pi}\circ i$ is {\it simple} if 
there exists a representative $\tilde{C}\in C(U,J^1({\mathbb R}^n,{\mathbb R}))$ of
 a extension of $i$ such that 
$\{ \tilde{C}_w| w\in  U\}$ is covered by 
finite orbits $[C_1],\ldots,[C_m]$ for some
 contact diffeomorphism germs $C_1,\ldots,C_m$ on $(J^1({\mathbb R}^n,{\mathbb R}),0)$. 

By an analogous way of Lemma \ref{simpleLagket:lem},
Lemma  \ref{simpleLag:lem} and Theorem \ref{genericclassLag:th} 
we have the followings:
\begin{lem}\label{simpleLegket:lem}
Let $\tilde{\pi} \circ i: (\tilde{\mathbb L},0)\rightarrow 
(J^1({\mathbb R}^n,{\mathbb R}),0)
\rightarrow  ({\mathbb R}^n\times{\mathbb R},0)$ 
be a stable reticular Legendrian map.
Then $\tilde{\pi} \circ i$ is simple if and only if there exist a 
neighbourhood $U_z$ of 
$z=j^{n+3}_0C(0)$ in $C^{n+3}(n)$ and $z_1,\ldots,z_m\in C^{n+3}(n)$ such 
that 
$U_z\subset [z_1]\cup\cdots\cup [z_m]$.
\end{lem}
\begin{lem}\label{simpleLeg:lem}
A stable reticular Legendrian  map $\tilde{\pi}\circ i$ is simple if and only if 
for a generating family 
$F(x,y,q,z)\in {\mathfrak M}(r;k+n+1)$ of $\tilde{\pi}\circ i$,
$F(x,y,0,0)\in {\mathfrak M}(r;k)^2$ is ${\cal K}$-simple singularity.
\end{lem}
\begin{tth}\label{genericclassLeg:th}
Let $r=0,n\leq 6$ or $r=1,n\leq 4$.
Let  $U$ be a neighbourhood of $0$ in 
$J^1({\mathbb R}^n,{\mathbb R})$.
Then there exists a residual set $O\subset C(U,J^1({\mathbb R}^n,{\mathbb R}))$
such that for any $\tilde{C}\in O$ and $w=(0,\ldots,0,p_{r+1},\ldots,p_n)\in U$,
the reticular Legendrian map $\tilde{\pi}\circ \tilde{C}_w|_{\tilde{\mathbb L}}$ 
is stable.
\end{tth}

In the case $r=0,n\leq 6$.  
A reticular 
Legendrian map $\tilde{\pi}\circ \tilde{C}_w|_{\tilde{\mathbb L}}$
for any $\tilde{C}\in O$ and $w \in U$ 
has a generating family which is a 
reticular ${\cal P}$-${\cal K}$-stable unfolding of one of 
$A_1,A_2,A_3,A_4,A_4,A_6,D^\pm_4,D_5,D^\pm_6,E_6$,
that is 
$F$ is stably reticular ${\cal P}$-${\cal K}$-equivalent to 
one of the following list:\\
$A_2:F(y_1,z)=y_1^2+z$,\\
$A_2:F(y_1,q_1,z)=y_1^3+q_1y_1+z$,\\
$A_3:F(y_1,q_1,q_2,z)= y_1^4+q_1y_1^2+q_2y_1+z$,\\
$A_4:F(y_1,q_1,q_2,q_3)= y_1^5+q_1y_1^3+q_2y_1^2+q_3y_1+z$,\\
$A_5:F(y_1,q_1,q_2,q_3,q_4)= 
 y_1^6+q_1y_1^4+q_2y_1^3+q_3y_1^2+q_4y_1+z$,\\
$A_6:F(y_1,q_1,q_2,q_3,q_4,q_5,z)= 
y_1^7+q_1y_1^5+q_2y_1^4+q_3y_1^3+q_4y_1^2+q_5y_1+z$,\\
$D^\pm_4:F(y_1,y_2,q_1,q_2,q_3,z)=
y_1^2y_2\pm y_2^3+q_1y_2^2+q_2y_2+q_3y_1+z$,\\
$D_5:F(y_1,y_2,q_1,q_2,q_3,q_4,z)=
y_1^2y_2+ y_2^4+q_1y_2^3+q_2y_2^2+q_3y_2+q_4y_1+z$,\\
$D^\pm_6:F(y_1,y_2,q_1,q_2,q_3,q_4,q_5,z)=
y_1^2y_2\pm  y_2^5+q_1y_2^4+q_2y_2^3+q_3y_2^2+q_4y_2
+q_5y_1+z$,\\
$E_6:F(y_1,y_2,q_1,q_2,q_3,q_4,q_5,z)=
y_1^3+ y_2^4+q_1y_1y_2^2+q_2y_1y_2+q_3y_2^2+q_4y_1
+q_5y_2+z$.\
\vspace{3mm}\\
In the case $r=1,\ n\leq 4$.   
A reticular Legendrian map $\tilde{\pi}\circ \tilde{C}_w|_{\tilde{\mathbb L}}$
for any $\tilde{C}\in O$ and $w \in U$
has a generating family which is a 
reticular ${\cal P}$-${\cal K}$-stable unfolding of one of 
$B_2,B_3,B_4,C^\pm_3,C_4,F_4$, that is 
$F$ is stably reticular ${\cal P}$-${\cal K}$-equivalent to 
one of the following list:\\
$B_2:F(x,q_1,z)= x^2+q_1x+z$,\\
$B_3:F(x,q_1,q_2,z)=x^3+q_1x^2+q_2x+z$,\\
$B_4:F(x,q_1,q_2,q_3,z)=x^4+q_1x^3+q_2x^2+q_1x+z$,\\
$C^\pm_3:F(x,y,q_1,q_2,z)=\pm xy+y^3+q_1y^2+q_2y+z$,\\
$C_4:F(x,y,q_1,q_2,q_3,z)=xy+y^4+q_1y^3+q_2y^2+q_3y+z$,\\
$F_4:F(x,y,q_1,q_2,q_3,z)=x^2+y^3+q_1xy+q_2x+q_3y+z$.\vspace{3mm}\\
\begin{figure}[htbp]
\[
\begin{array}{cc}
\includegraphics*[width=7cm,height=7cm]{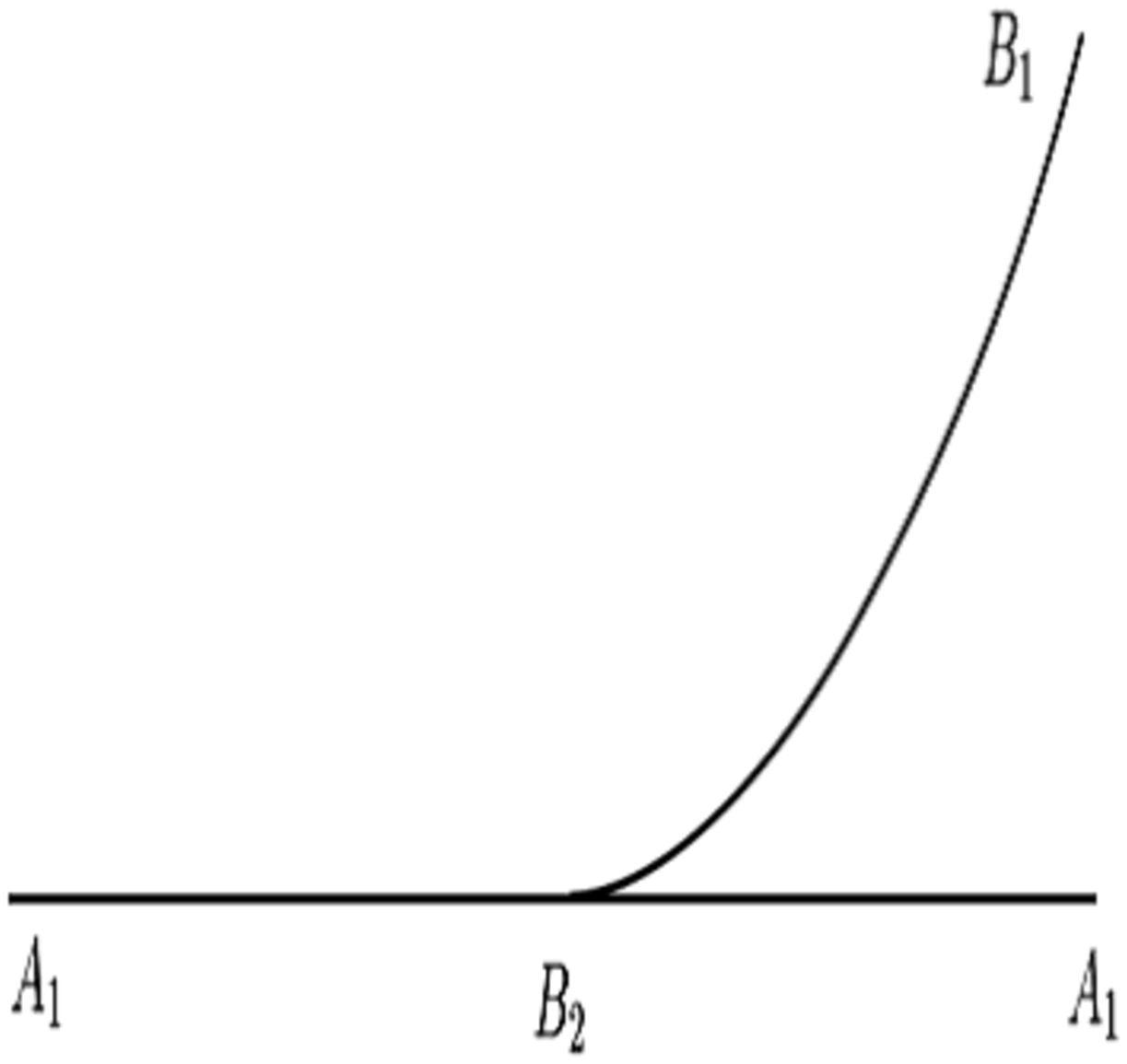}&
\includegraphics*[width=7cm,height=7cm]{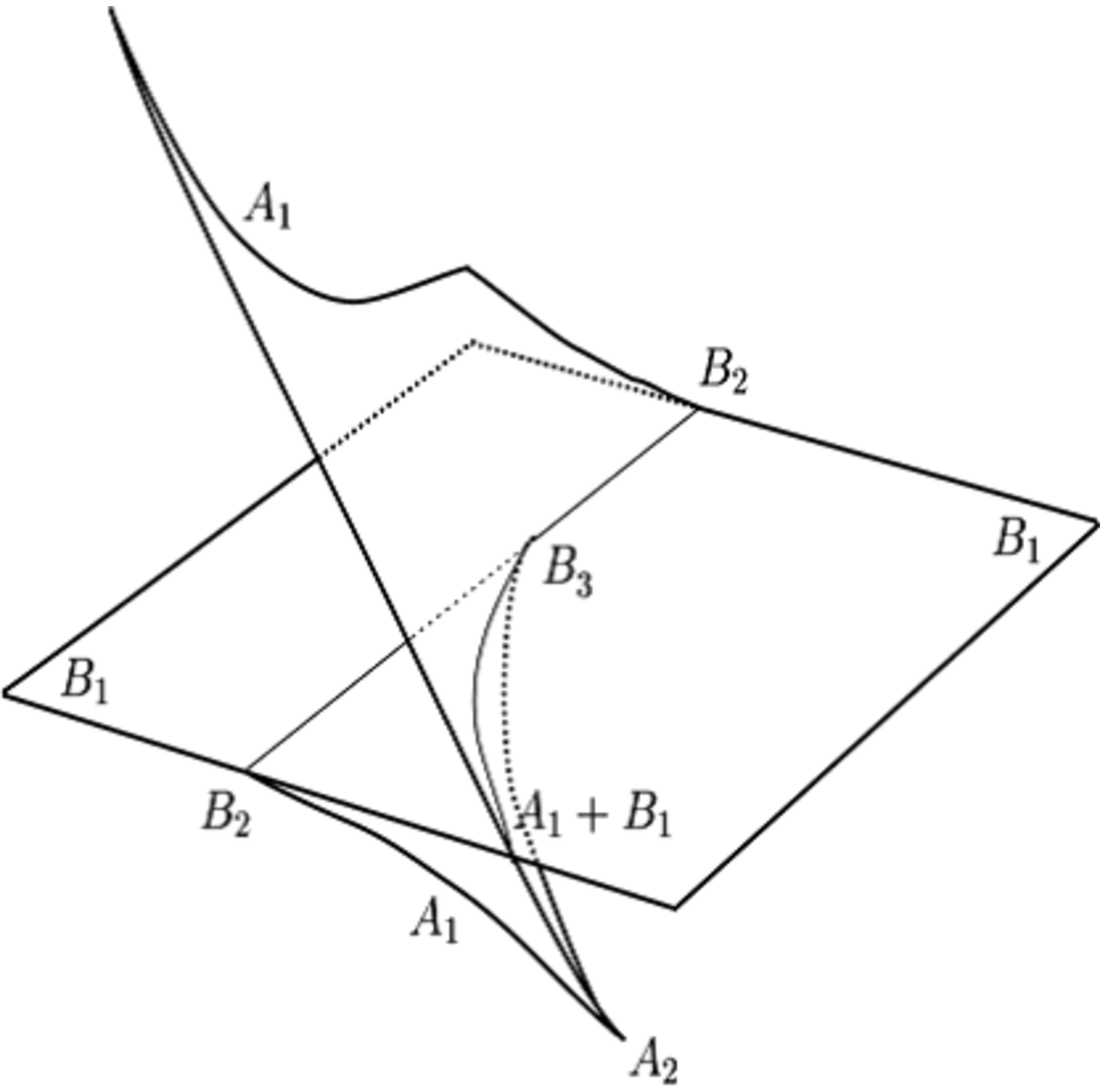}
\end{array}\]
\caption{the wavefronts $B_2$ and $B_3$}
\[
\begin{array}{cc}
\includegraphics*[width=7cm,height=7cm]{FC+_3.eps}
&
\includegraphics*[width=7cm,height=7cm]{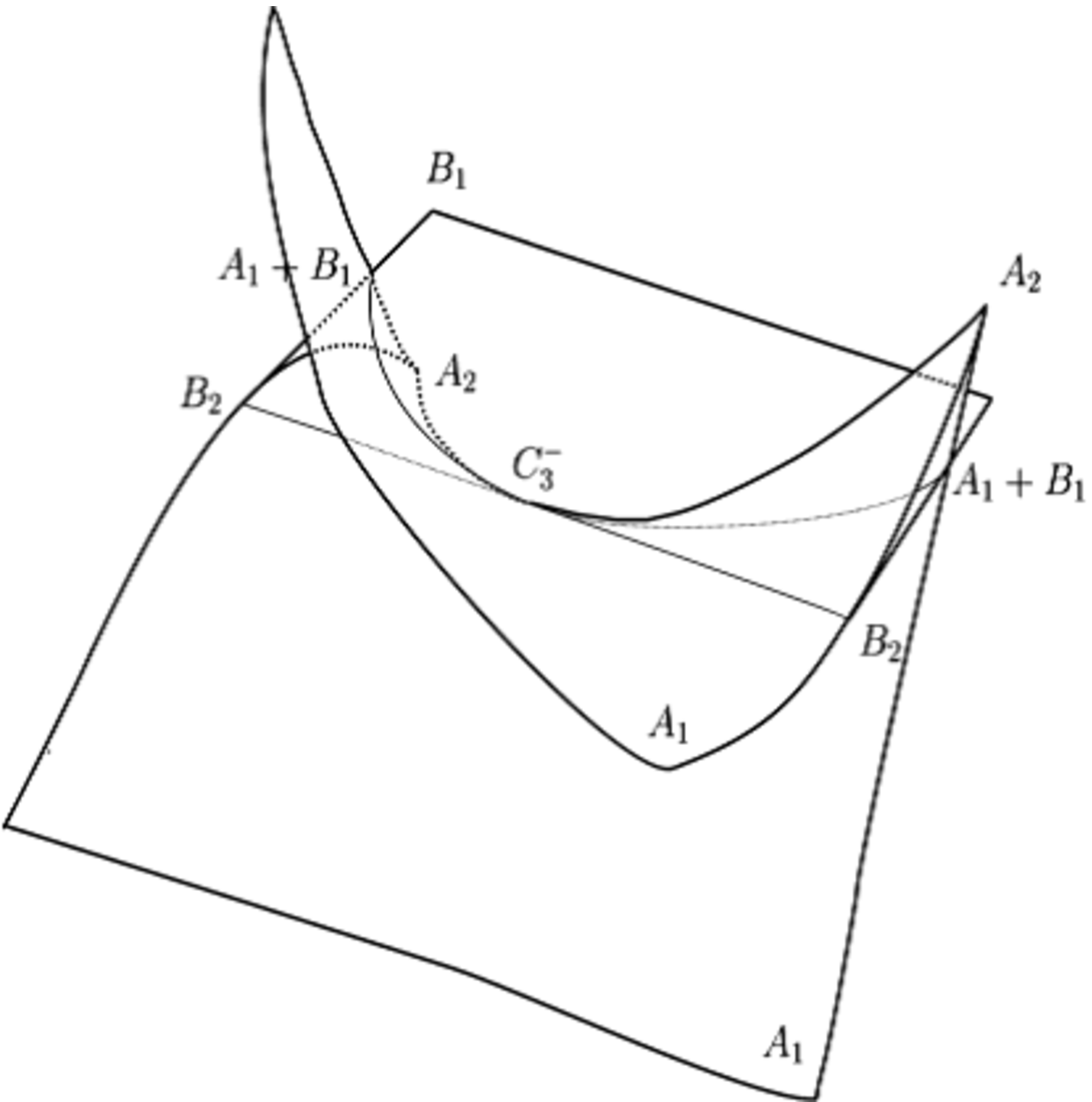}
\end{array}\]
\caption{the wavefronts $C^\pm_3$}
\end{figure}%


\begin{thebibliography}{99}
\bibitem{arnold:text}
V.I.Arnold, S.M. Gusein-Zade \& A.N. Varchenko, {\it Singularities of differential maps,} Vol.I, Birkh user, Basel, 1986.
\bibitem{golubitsky:text}
M.Golubitsky, V.Guillemin, {\it Stable Mappings and Their Singularities}, 
Graduate Text in Math. 14, Springer-Verlag, New York 1973. 
\bibitem{izumiya:text}
S.Izumiya G.Ishikawa, {\it Singularity Theory and their applications}, Kyouritsu shuppan, 1998 (in Japanese)
\bibitem{izumiya1}
S.Izumiya, {\it Semi-local Classification of Geometric Singularities
for Hamilton-Jacobi Equations},  J. Differentials Geom., {\bf 118} (1995), 
pp.166-193
\bibitem{Broc:text}
Th:br\"ocker, {\it Differentiable Germs and Catastrophes},
London Math. Soc. Lecture Note ser., 17 Cambridge Univ. Press, 1975.
\bibitem{janeszko3}
S.Janeszko, {\it Generalized Luneburg canonical varieties and vector fields on 
quasicaustics}, J. Math. Phys. 31:4(1990), pp.997-1009.
\bibitem{janich1}
K.J\"anich, {\it Caustics and catastrophes}, Math. Ann., 209(1974), pp.161-180.\bibitem{gakui}
T.Tsukada, {\it Reticular Lagrangian, legendrian Singularities and
their applications}, PhD Thesis, Hokkaido University, 1999.
\bibitem{retLag}
T.Tsukada, {\it Reticular Lagrangian Singularities}, 
The Asian J. of Math., vol 1, {\bf 3} (1997), pp.572-622.
\bibitem{retLeg}
T.Tsukada, {\it Reticular Legendrian Singularities}, 
The Asian J. of Math., vol 5, {\bf 1} (2001), pp.109-127.
\bibitem{ftsing}
G.Wassermann, {\it Stability of unfolding}, 
Lecture note in mathematics, 393.
\bibitem{spsing}
G.Wassermann, {\it Stability of unfolding in space and time},
Acta Math., {\bf 135} (1975), pp.57-128.
\bibitem{zaka1}
V.M.Zakalyukin, {\it Reconstruction of fronts and caustics depending on 
a parameter and versality of mappings}, J. of Soviet Math., {\bf 27} (1984), 
pp.2713-2735
\end{thebibliography}
\end{document}